\documentclass[12pt, leqno]{amsart}

\usepackage{lscape}
\usepackage{graphicx}
\usepackage{amssymb,amsmath,amsthm,amscd}
\usepackage{mathrsfs}
\usepackage{enumerate}
\usepackage[usenames,dvipsnames]{color}
\usepackage[colorlinks=true, pdfstartview=FitV,
 linkcolor=blue,citecolor=blue,urlcolor=blue]{hyperref}

\usepackage[all]{xy}

\allowdisplaybreaks[4]
\usepackage{oldgerm}

\newcommand{\nc}{\newcommand}
\numberwithin{equation}{section}

\newenvironment{jaune}{\relax\color{Orchid}}{\hspace*{.5ex}\relax}

\newcommand{\bj}{\begin{jaune}}
\newcommand{\ej}{\end{jaune}}

\theoremstyle{plain}
\newtheorem{lemma}{Lemma}[section]
\newtheorem{prop}[lemma]{Proposition}
\newtheorem{theorem}[lemma]{Theorem}
\newtheorem*{maintheorem}{Main Theorem}
\newcommand{\Prop}{\begin{prop}}
\newcommand{\enprop}{\end{prop}}
\newcommand{\Lemma}{\begin{lemma}}
\newcommand{\enlemma}{\end{lemma}}
\newcommand{\Th}{\begin{theorem}}
\newcommand{\enth}{\end{theorem}}
\newtheorem{corollary}[lemma]{Corollary}
\newcommand{\Cor}{\begin{corollary}}
\newcommand{\encor}{\end{corollary}}
\newtheorem{definition}[lemma]{Definition}
\newtheorem*{conjecture}{Conjecture}
\newcommand{\Def}{\begin{definition}}
\newcommand{\edf}{\end{definition}}
\newtheorem{sublemma}[lemma]{Sublemma}
\newcommand{\Sublemma}{\begin{sublemma}}
\newcommand{\ensub}{\end{sublemma}}

\theoremstyle{definition}
\newtheorem{remark}[lemma]{Remark}

\newtheorem{Convention}[lemma]{Convention}
\newcommand{\Conv}{\begin{Convention}}
\newcommand{\enconv}{\end{Convention}}
\nc{\Con}{\begin{conjecture}}
\nc{\encon}{\end{conjecture}}
\nc{\Rem}{\begin{remark}}
\nc{\enrem}{\end{remark}}
\newcommand{\C}{{\mathbb C}}
\newcommand{\Q}{\mathbb {Q}}

\newcommand{\Z}{{\mathbb Z}}
\newcommand{\B}{{\mathbf{B}}}
\newcommand{\A}{{\mathbf A}}

\newcommand{\Db}[1][G]{{\operatorname{D^b_{#1}}}}

\newcommand{\R}{{\rm R}}

\newcommand{\one}{{\bf{1}}}
\newcommand{\seteq}{\mathbin{:=}}

\newcommand{\hd}{{\operatorname{hd}}}
\newcommand{\soc}{{\operatorname{soc}}}

\newcommand{\g}{{\mathfrak{g}}}

\newcommand{\Hom}{\operatorname{Hom}}

\newcommand{\isoto}[1][]{\mathop{\xrightarrow%
[{\raisebox{.3ex}[0ex][.3ex]{$\scriptstyle{#1}$}}]%
{{\raisebox{-.6ex}[0ex][-.6ex]{$\mspace{2mu}\sim\mspace{2mu}$}}}}}
\newcommand{\tensor}{\otimes}

\newcommand{\eq}{\begin{eqnarray}}
\newcommand{\eneq}{\end{eqnarray}}

\newcommand{\hs}{\hspace*}

\newcommand{\To}[1][{\hs{2ex}}]{\xrightarrow{\,#1\,}}

\newcommand{\eqn}{\begin{eqnarray*}}
\newcommand{\eneqn}{\end{eqnarray*}}
\newcommand{\on}{\operatorname}
\newcommand{\Ker}{\on{Ker}}

\newcommand{\bna}{\be[{\rm(a)}]}

\newcommand{\QED}{\end{proof}}
\newcommand{\Proof}{\begin{proof}}

\newcommand{\soplus}{\mathop{\mbox{\normalsize$\bigoplus$}}\limits}
\newcommand{\sodot}{\mathop{\mbox{\normalsize$\bigodot$}}\limits}
\newcommand{\sotimes}{\mathop{\mbox{\normalsize$\bigotimes$}}\limits}

\newcommand{\id}{\on{id}}
\newcommand{\ba}{\begin{array}}
\newcommand{\ea}{\end{array}}

\newcommand{\monoto}{\rightarrowtail}

\newcommand{\set}[2]{\left\{#1 \mid #2 \right\}}

\newcommand{\eqsub}{\begin{subequations}\begin{eqnarray}}
\newcommand{\eneqsub}{\end{eqnarray}\end{subequations}}

\newcommand{\ol}{\overline}

\nc{\la}{\lambda}
\nc{\lam}{\lambda}
\nc{\U}[1][\g]{U_q(#1)}
\nc{\te}{\tilde{e}}
\nc{\tei}{\tilde{e}_i}
\nc{\tf}{\tilde{f}}
\nc{\tfi}{\tilde{f}_i}
\nc{\tU}{\widetilde U_q(\g)}
\nc{\tE}{\tilde{E}}
\nc{\tF}{\widetilde{\F}}
\nc{\tK}{\widetilde{K}}

\nc{\tk}{\tilde{k}}
\nc{\tkone}{\tk_{\ol{1}}}
\nc{\teone}{\tilde{e}_{\ol{1}}}
\nc{\tfone}{\tilde{f}_{\ol{1}}}

\nc{\teibar}{\tilde{e}_{\ol{i}}} \nc{\tfibar}{\tilde{f}_{\ol{i}}}
\nc{\tki}{{\tk}_{\ol {i}}}

\nc{\BZ}{{\mathbb{Z}}}
\nc{\al}{\alpha}
\nc{\qs}{{q}}
\nc{\lan}{\langle}
\nc{\ran}{\rangle}
\nc{\re}{{\mathrm{re}}}
\nc{\wt}{\operatorname{wt}}
\nc{\ch}{\operatorname{ch}}
\nc{\Um}[1][\g]{U^-_q(#1)}
\nc{\Ue}{U^+_q(\g)}
\nc{\eps}{\varepsilon}
\nc{\vphi}{\varphi}
\nc{\sphi}{\varphi^*}
\nc{\seps}{\varepsilon^*}

\nc{\nn}{\nonumber}
\def\max{{\mathop{\mathrm{max}}}}
\nc{\vp}{\varpi}
\nc{\cls}{{\operatorname{cl}}}
\nc{\Wt}{{\operatorname{Wt}}}
\nc{\Us}{U'_q(\g)}
\nc{\La}{\Lambda}
\nc{\tLa}{\widetilde\Lambda}
\nc{\ro}{{\rm(}}
\nc{\rf}{{\rm)}}
\nc{\norm}{{\mathrm{norm}}}
\nc{\qbox}{\quad\mbox}
\nc{\braid}{{\mathfrak{B}}}
\nc{\Ad}{\operatorname{Ad}}
\nc{\Aut}{\operatorname{Aut}}
\nc{\dt}[1]{\tilde{\tilde #1}}
\nc{\Sn}{S^{{\mathrm{norm}}}}
\nc{\aff}{{\rm{aff}}}
\nc{\rk}{{\mathrm{rk}}}
\nc{\tP}{\widetilde{P}}
\nc{\tW}{\widetilde{W}}
\nc{\Dyn}{\mathrm{Dyn}}
\nc{\tD}{\widetilde{\Delta}}
\nc{\height}[1]{{\operatorname{ht}}(#1)}
\nc{\bl}{\bigl(}
\nc{\br}{\bigr)}
\nc{\Hecke}{\mathrm{H}}
\nc{\HA}{\Hecke^{\mathrm{A}}}
\nc{\HB}{\Hecke^{\mathrm{B}}}
\newcommand{\scbul}{{\,\raise1pt\hbox{$\scriptscriptstyle\bullet$}\,}}
\nc{\vac}{{\phi}}
\nc{\Bt}{\B_\theta(\g)}
\nc{\be}{\begin{enumerate}}
\nc{\ee}{\end{enumerate}}
\nc{\low}{{\mathrm{low}}}
\nc{\upper}{{\mathrm{up}}}
\nc{\Zodd}{\Z_{\mathrm{odd}}}
\nc{\Ft}[1][n]{\mathbb{P}\mathrm{ol}_{#1}}
\nc{\Ftf}[1][n]{\widetilde{\mathbb{P}\mathrm{ol}}_{#1}}
\nc{\KA}{\on{K}^{\mathrm{A}}}
\nc{\KB}{\on{K}^{\mathrm{B}}}
\nc{\Res}{\on{Res}}
\nc{\Fc}[1][{n,m}]{\mathbf{F}_{#1}}
\nc{\tphi}{\tilde{\varphi}}
\nc{\CO}{\mathscr{O}}
\nc{\inte}{\mathrm{int}}
\nc{\Oint}{\mathcal{O}^{\ge0}_{\inte}}
\nc{\vs}{\vspace*}
\nc{\tLt}{\widetilde{L}}
\nc{\tL}{\widetilde{\Lambda}}
\nc{\tu}{\tilde{u}}
\nc{\noi}{\noindent}
\nc{\heigh}{\mathfrak{t}}
\nc{\lowest}{\mathfrak{l}}
\nc{\rootl}{\mathsf{Q}}
\nc{\cl}{\colon}
\nc{\uqpg}{U'_q(\mathfrak g)}
\nc{\uq}{\uqpg}
\nc{\Oh}{\widehat{\mathcal{O}}}

\nc{\KLR}{KLR algebra}
\nc{\KLRs}{KLR algebras}
\nc{\cor}{\mathbf{k}}
\nc{\cora}{{\cor(A)}}
\nc{\haut}{\mathrm{ht}}
\nc{\tens}{\mathop\otimes}
\nc{\gmod}{\mbox{-$\mathrm{gmod}$}}
\nc{\gMod}{\mbox{-$\mathrm{gMod}$}}
\nc{\proj}{\mbox{-$\mathrm{proj}$}}
\nc{\gproj}{\mbox{-$\mathrm{gproj}$}}
\nc{\smod}{\mbox{-$\mathrm{mod}$}}
\nc{\Mod}{\mbox{-$\mathrm{Mod}$}}
\nc{\h}{\mathfrak h}
\nc{\Rnorm}{R^{\rm{norm}}}

\nc{\Vhat}{\widehat{V}}
\nc{\F}{\mathcal{F}}

\def\T{{\mathcal T}}

\nc{\fd}[1][A]{\on{\mathrm{flat.dim}_{#1}}}
\nc{\bP}{{\mathbb{P}}}
\nc{\bPh}{\widehat{\mathbb{P}}}
\nc{\bK}[1][{n}]{\widehat{\mathbb{K}}_{#1}}
\nc{\bV}[1][{n}]{\widehat{V}^{\otimes{#1}}}
\nc{\bVK}[1][{n}]{\widehat{V}^{\otimes{#1}}_{\widehat{\mathbb{K}}}}
\nc{\hV}{\widehat{V}}
\nc{\opp}{\mathrm{opp}}
\nc{\col}{\colon}
\nc{\bnum}{\be[{\rm(i)}]}
\nc{\oep}{\epsilon}
\nc{\qtext}{\quad\text}
\nc{\qtextq}[1]{\quad\text{#1}\quad}
\nc{\longtwoheadrightarrow}[1][]{\xymatrix{\ar@{->>}[r]^-{{#1}}&}}
\nc{\epiTo}[1][]{\longtwoheadrightarrow[{#1}]}
\nc{\epito}{\twoheadrightarrow}
\nc{\monoTo}[1][]{\xymatrix{\ar@{>->}[r]^-{{#1}}&}}
\nc{\sym}{\mathfrak{S}}
\nc{\inp}[1]{{({#1})_{\mathrm{n}}}}
\nc{\rtl}{\rootl}
\nc{\wtd}{\widetilde}
\nc{\etens}{\boxtimes}
\nc{\ds}[1]{\mathrm{d}(#1)}
\nc{\rmat}[1]{{\mathbf{r}}_%
{\mspace{-2mu}\raisebox{-.6ex}{${\scriptstyle{#1}}$}}}
\nc{\rmats}[1]{{\mathbf{r}}_%
{\mspace{-2mu}\raisebox{-.6ex}{${\scriptscriptstyle{#1}}$}}}
\nc{\shc}{\mathcal{C}}
\nc{\shs}{\mathcal{S}}
\nc{\Fct}{{\on{Fct}}}
\nc{\tC}{\widetilde{\shc}}
\nc{\Zp}{\Z_{\ge0}}
\nc{\tPhi}{\widetilde{\Phi}}
\nc{\tT}{{\widetilde{\T}}}
\nc{\Ob}{\on{Ob}}
\nc{\bwr}{\mbox{\large$\wr$}}
\nc{\Img}{\on{Im}}
\nc{\Ab}{\mathcal{A}^{\mathrm{big}}}
\nc{\Sb}{\mathcal{S}^{\mathrm{big}}}
\nc{\As}{\mathcal{A}}
\nc{\Ss}{\mathcal{S}}
\nc{\ntens}{\widetilde{\otimes}}
\nc{\hR}{\widehat{R}}
\nc{\nconv}{\mathop{\mbox{\large $\odot$}}}
\nc{\snconv}{\mbox{\scriptsize$\odot$}}
\nc{\ts}{\tilde{s}}
\nc{\sho}{\mathcal{O}}
\nc{\bc}{\begin{cases}}
\nc{\ec}{\end{cases}}
\nc{\slnh}{{\widehat{\mathfrak{sl}}_N}}
\nc{\UA}{U_q'(\slnh)}
\nc{\KR}{R_K}
\nc{\cQ}{\mathcal{Q}}
\nc{\Irr}{\mathcal{I}rr}
\nc{\tQ}{\widetilde{\cQ}}
\nc{\bs}{\mathbf{s}}
\nc{\bL}{\mathbb{L}}
\nc{\tg}{\tilde{g}}

\nc{\conv}{\mathbin{\mbox{\large $\circ$}}}
\nc{\shconv}{\mathbin{\large\diamond}}
\nc{\hconv}{\mathbin{\mbox{\Large $\shconv$}}}

\nc{\Rm}{R^{\mathrm{ren}}}

\nc{\bQ}{\ol{Q}}
\renewcommand{\Im}{\on{Im}}

\nc{\de}{\on{\textfrak{d}}}

\nc{\xmono}{\ar@{>->}}
\nc{\xepi}{\ar@{->>}}
\nc{\db}[1]{\raisebox{-.5ex}[2ex][1.8ex]{$#1$}}
\nc{\wb}[1]{\mbox{$\rule[-1.1ex]{0ex}{2ex}#1$}}
\nc{\univ}{\mathrm{univ}}
\nc{\rM}{{}^*\mspace{-2mu}M}
\nc{\lM}{M^*}
\nc{\uqm}{\uq\smod}
\nc{\tR}{\widetilde{R}_{\gamma,\beta}}
\nc{\tx}{\tilde{x}}
\nc{\bi}{\mathbf{i}}
\nc{\ttau}{\widetilde{\tau}}

\nc{\tEnd}{\on{\widetilde{E}nd}}
\nc{\tHom}{\on{\widetilde{H}om}}

\nc{\K}{{J}}
\nc{\Kex}{{\K}_{\mathrm{ex}}}
\nc{\Kfr}{{\K}_{\mathrm{f\mspace{.01mu}r}}}
\nc{\coro}{\cor}
\nc{\tB}{\widetilde{B}}
\nc{\seed}{\mathscr{S}}

\nc{\up}{\mathrm{up}}
\nc{\bfa}{\mathbf{a}}
\newcommand{\wB}{\widetilde{B}}

\newlength{\mylength}
\title
{Monoidal categorification of cluster algebras}

\author[S.-J. Kang, M. Kashiwara, M. Kim, Se-jin Oh]{Seok-Jin Kang$^{1}$, Masaki Kashiwara$^{2}$, Myungho Kim, Se-jin Oh$^3$}

\address{Department of Mathematical Sciences
         and
         Research Institute of Mathematics \\
         Seoul National University \\ Seoul 151-747, Korea}
         \email{sjkang@math.snu.ac.kr}

\address{Research Institute for Mathematical Sciences \\
          Kyoto University \\ Kyoto 606-8502, Japan \\
          \& Department of Mathematical Sciences
         and
         Research Institute of Mathematics \\
         Seoul National University \\ Seoul 151-747, Korea}
         \email{masaki@kurims.kyoto-u.ac.jp}

\address{School of Mathematics, Korea Institute for Advanced Study \\ Seoul 130-722, Korea}
         \email{mhkim@kias.re.kr}

\address{Department of Mathematical Sciences
         and
         Research Institute of Mathematics \\
         Seoul National University \\ Seoul 151-747, Korea}
         \email{sj092@snu.ac.kr}

\thanks{$^1$ This work was supported by NRF grant \# 2014021261
and NRF grant \# 2013055408
}
\thanks{$^2$ This work was supported by Grant-in-Aid for
Scientific Research (B) 22340005, Japan Society for the Promotion of
Science.}
\thanks{$^3$ This work was supported by BK21 PLUS SNU Mathematical Sciences Division}


\date{December 28, 2014}

\begin{document}

\begin{abstract}

We give a definition of monoidal categorifications of quantum
cluster algebras and provide a criterion for a monoidal category of
finite-dimensional graded $R$-modules to become a monoidal
categorification of a quantum cluster algebra,
where $R$ is a symmetric Khovanov-Lauda-Rouquier algebra. 
Roughly speaking, this criterion asserts that
a quantum monoidal seed can be mutated successively in all the directions
once the first-step mutations are possible.
In the course of the study, we also give a proof of
a conjecture of Leclerc on
the product of upper global basis elements.

\end{abstract}
\maketitle

\section*{Introduction}
The purpose of this paper is to give a definition of
a monoidal categorification of a quantum cluster algebra
and to provide a criterion for a monoidal category
to be a monoidal categorification.

The notion of cluster algebras was introduced by Fomin and Zelevinsky in
\cite{FZ02} for studying total positivity and upper global bases.
Since their introduction, connections and applications have been
discovered in various fields of mathematics including 
representation theory, Teichm\"uller theory, tropical geometry, integrable
systems, and Poisson geometry.

A cluster algebra is a $\Z$-subalgebra of a rational function field
given by a set of generators,
called the {\it cluster variables}. These generators are grouped
into overlapping subsets, called {\it clusters}, and the clusters are defined
inductively by a procedure called {\it mutation} from the {\it
initial cluster} $\{ X_i\}_{1 \le i \le r}$, which is controlled by
an exchange matrix $\wB$. We call a monomial of cluster
variables in one cluster {\it a cluster monomial}.

Fomin and Zelevinsky proved that every cluster variable is a Laurent
polynomial of the initial cluster $\{ X_i\}_{1 \le i \le r}$ and
they conjectured that this Laurent polynomial has positive
coefficients (\cite{FZ02}). This {\it positivity conjecture} was
proved by Lee and Schiffler in the {\it skew-symmetric} cluster algebra case
in \cite{LS13}.
The {\it linearly independence
conjecture} on cluster monomials was proved 
in the skew-symmetric cluster algebra case
in \cite{CKLP12}.

The notion of quantum cluster algebras, introduced by Berenstein and Zelevinsky in
\cite{BZ05}, can be considered as a $q$-analogue of cluster algebras.
The commutation relation among the cluster variables is determined by
a skew-symmetric matrix $L$. 
As in the cluster algebra case, every cluster variable belongs to
$\Z[q^{\pm 1/2}][X_i^{\pm 1}]_{1 \le i \le r}$ (\cite{BZ05}), and is
expected to be an element of $\Z_{\ge0}[q^{\pm 1/2}][X_i^{\pm 1}]_{1 \le i \le r}$, 
which is referred to as the {\it quantum
positivity conjecture} (cf.\ \cite[Conjecture 4.7]{DMSS}).
In \cite{KF14}, Kimura and Qin proved the quantum positivity conjecture for quantum cluster algebras containing {\it acyclic} seed
and specific coefficients.

In a series of papers \cite{GLS11,GLS07,GLS}, Gei\ss, Leclerc and
Schr{\"o}er showed that the quantum unipotent coordinate algebra
$A_q(\mathfrak{n}(w))$, associated with a symmetric quantum group $U_q(\g)$
and its Weyl group element $w$, has a 
skew-symmetric quantum cluster algebra structure whose initial cluster
consists of {\it quantum minors}. In \cite{Kimu12}, Kimura
proved that $A_q(\mathfrak{n}(w))$ is {\it compatible} with the
upper global basis $\B^{\upper}$ of $A_q(\mathfrak{n})$; i.e., the
set $\B^{\upper}(w) \seteq A_q(\mathfrak{n}(w)) \cap \B^{\upper}$
is a basis of $A_q(\mathfrak{n}(w))$. Thus, with a
result of \cite{CKLP12}, one can expect that every cluster monomial
of $A_q(\mathfrak{n}(w))$ is contained in the upper global basis
$\B^{\upper}(w)$, which is named {\it the quantization conjecture} by Kimura
(\cite{Kimu12}).

In \cite{HL10}, Hernandez and Leclerc introduced the notion of {\it
a monoidal categorification of  a cluster algebra}. 
We say that a simple object $S$ of a monoidal category
$\shc$ is {\it real} if $S \tens S$ is simple. We say that
a simple object $S$ is {\it prime} if there exists no non-trivial
factorization $S \simeq S_1 \tens S_2$. They say that $\shc$ is
a monoidal categorification of a cluster algebra $A$ if
the Grothendieck ring of $\shc$ is isomorphic to
$A$ and if 

\vs{1.5ex}
\hs{0ex}\parbox{80ex}{

\begin{enumerate}
\setlength{\itemsep}{3pt}
\item[{\rm (M1)}] the  cluster monomials of $A$ are the classes of real simple objects of $\shc$,
\item[{\rm (M2)}] the cluster variables of $A$ are the classes of real simple prime objects of $\shc$.
\end{enumerate}}

\vs{1.5ex}
\noi
(Note that the above version is
weaker than the original definition of the monoidal categorification in
\cite{HL10}.) They proved that certain categories 
of modules over symmetric quantum affine algebras
$U_q'(\g)$ give monoidal categorifications of cluster algebras.
Nakajima extended it to the cases of the cluster algebras of type $A,D,E$
(\cite{Nak11}) (see also \cite{HL13}).

Once a cluster algebra $A$ has a monoidal categorification,
the positivity of cluster
variables of $A$ and the linear independency of cluster monomials of $A$ follow
(see \cite[Proposition 2.2]{HL10}).

In this paper, we will refine their notion of
monoidal categorifications including the quantum cluster algebra case.

\medskip

The Khovanov-Lauda-Rouquier (or simply KLR) algebras, introduced by
Khovanov-Lauda \cite{KL09,KL11} and Rouquier \cite{R08}
independently, are a family of $\Z$-graded algebras which
categorifies the negative half $U_q^-(\g)$ of a {\it symmetrizable}
quantum group $U_q(\g)$. More precisely, there exists a family of
algebras $\{ R(-\beta) \}_{\beta \in \rtl^-}$ such that the Grothendieck
ring of $R \gmod \seteq \bigoplus_{\beta \in \rtl^-}R(-\beta)\gmod$, the direct sum
of the categories of finite-dimensional graded $R(-\beta)$-modules, is
isomorphic to the integral form $A_q(\mathfrak{n})_{\Z[q^{\pm1}]}$ of
$A_q(\mathfrak{n}) \simeq U_q^-(\g)$. Here the tensor functor $\tens$
of the monoidal category $R \gmod$ is given by
the convolution product $\conv$, and the action of $q$ is given by
the grading shift functor. In \cite{VV09, R11},
Varagnolo-Vasserot and Rouquier
proved that the upper global basis $\B^\upper$ of $A_q(\mathfrak{n})$
corresponds to the
set of the classes of all {\it self-dual} simple modules of $R
\gmod$ under the assumption that $R$ is associated with a {\it
symmetric} quantum group $U_q(\g)$.

Combining works of \cite{GLS,Kimu12,VV09}, the quantum unipotent
coordinate algebra $A_q(\mathfrak{n}(w))$ associated with a
symmetric quantum group $U_q(\g)$ and a Weyl group element $w$
is isomorphic to the Grothendieck group of
the monoidal abelian full subcategory $\shc_w$ of $R \gmod$
satisfying the following properties: {\rm (i)} $\shc_w$ is
stable under extensions and grading shift functor, {\rm (ii)} the
composition factors of $M \in \shc_w$ are contained
in $\B^{\upper}(w)$. However it is not evident the conditions (M1) and (M2)
are satisfied. The purpose of this paper is provide a theoretical background 
in order to prove (M1) and (M2). In the forthcoming paper, we 
prove (M1) and (M2) as an application of the results of the present paper.

\medskip

In this paper, we first continue the work of \cite{KKKO14} about the
convolution products, heads and socles of graded modules over symmetric KLR algebras. 
One of the main results in \cite{KKKO14} is that the convolution
product $M \conv N$ of a real simple $R(\beta)$-module $M$ and a
simple $R(\gamma)$-modules $N$ has a unique simple quotient and a unique
simple submodule. Moreover, if $M \conv N \simeq N \conv M$ up to a
grading shift, $M \conv N$ is simple. 
In such a case we say that $M$ and $N$ {\it
commute}. The main tool of \cite{KKKO14} was an R-matrix
$\rmat{M,N}$, constructed in \cite{K^3}, which is a homogeneous
homomorphism from $M \conv N$ to $N \conv M$ of degree
$\La(M,N)$. 
We define the integers
\begin{align*}
\tLa(M,N) \seteq \frac{1}{2}\bl\La(M,N)+(\beta,\gamma)\br, \quad
\de(M,N)\seteq \frac{1}{2}\bl\La(M,N)+\La(N,M)\br
\end{align*}
and study the representation theoretic meaning of the integers $\La(M,N)$, $\tLa(M,N)$ and $\de(M,N)$.

We then prove Leclerc's first conjecture (\cite{L03}) 
on the multiplicative
structure of elements in $\B^\upper$ when the generalized Cartan matrix is
 symmetric (Theorem \ref{th:leclerc} and
Theorem \ref{th:head}). Theorem \ref{th:head} is due to McNamara 
(\cite[Lemma 7.5]{Mc14}) and the authors thank him for informing us his result.

 We say that $b \in \B^\upper$ is {\it real} if
$b^2 \in q^\Z\,\B^\upper\seteq\bigsqcup_{n\in\Z}q^n\B^\upper$.
\begin{conjecture}[{\cite[Conjecture 1]{L03}}] Let $b_1$ and $b_2$ be elements in $\B^\upper$ such that $b_1$ is real and $b_1b_2 \not\in q^\Z\B^\upper$.
Then the expansion of $b_1b_2$ with respect to $\B^\upper$ is of the
form
$$ b_1b_2= q^m b' + q^s b'' + \sum_{c \ne b',b''}\gamma^c_{b_1,b_2}(q)c,$$
where $b' \ne b''$, $m,s \in \Z$, $m < s$, $\gamma^c_{b_1,b_2}(q)\in
\Z[q^{\pm 1}]$, and for all $c \in \B^\upper$ such that
$\gamma^c_{b_1,b_2}(q) \ne 0$
$$  \gamma^c_{b_1,b_2}(q) \in q^{m+1}\Z[q] \cap q^{s-1}\Z[q^{-1}].$$
\end{conjecture}

More precisely, we prove that $q^mb'$ and $q^sb''$ correspond to the
simple head 
and the simple socle of $M\conv N$, respectively, 
when $b_1$ corresponds to a real
simple module $M$ and $b_2$ corresponds to a simple module $N$.

In the second part of this paper, we provide an algebraic
framework for monoidal categorifications of cluster algebras and
quantum cluster algebras. Let us focus on the quantum cluster algebra
case in this introduction.

Let $\K$ be a finite index set with a decomposition $\Kex
\sqcup \Kfr$ and let $\rtl$ be a free abelian group with a
symmetric bilinear form $( \ , \ )$ such that $(\beta,\beta) \in
2\Z$ for all $\beta \in \rtl$. Let $\shc$ be an abelian monoidal
category with a {\it grading shift functor} $q$ 
(see \eqref{cond:quantum monoidal category}) and
a decomposition $\shc= \bigoplus_{\beta \in \rtl}
\shc_\beta$.
The quadruple $\seed \seteq
(\{ M_i\}_{i \in \K}, L,\widetilde{B},D)$ is called a {\it quantum
monoidal seed}\/ when it consists of {\rm (i)} an integer-valued $\K
\times \K_{{\rm ex}}$-matrix $\widetilde{B}$ with the skew-symmetric
principal part, {\rm (ii)} an integer-valued skew-symmetric $\K
\times \K$-matrix $L=(\la_{ij})_{i,j\in\K}$ such that $(L,\wB)$ is {\it compatible} with
the integer $2$, {\rm (iii)} the set of real simple modules $\{ M_i \}_{i
\in \K}$ in $\shc$ such that $M_i \tens M_j \simeq
q^{\lambda_{ij}}M_j M_i$ and $M_{i_1} \tens \cdots \tens M_{i_\ell}$
is simple for any $i,j \in \K$ and $(i_1,\ldots,i_\ell) \in \K^\ell$,
{\rm (iv)} $D=({\rm wt}(M_i))_{i \in \K} \in \rtl^\K$, and it satisfies certain
conditions (Definition \ref{def:quantum monoidal seed}).
Here we set $\wt(M)=\beta$ for $M \in \shc_\beta$.

For $k \in \K_{{\rm ex}}$, we say that $\seed$ {\it admits a mutation
in direction $k$} if there exists a simple object $M_k'$ such that
there exist exact sequences in $\shc$
\begin{align*}
&0 \To q \sodot_{b_{ik} >0} M_i^{\snconv b_{ik}} \To q^{m_k} M_k
\tensor M_k' \To
 \sodot_{b_{ik} <0} M_i^{\snconv (-b_{ik})} \To 0, \\
&0 \To q \sodot_{b_{ik} <0} M_i^{\snconv(-b_{ik})} \To q^{m_k'} M_k'
\tensor M_k \To
  \sodot_{b_{ik} >0} M_i^{\snconv b_{ik}} \To 0,
\end{align*}
and its {\it mutation} $\mu_k(\seed) \seteq (\{ \mu_k(M)_i\}_{i \in
\K}, \mu_k(L),\mu_k(\widetilde{B}),\mu_k(D))$ is again a quantum
monoidal seed. Here $\sodot$ denotes a tensor product with some
grading shift (see \eqref{eq:balpro}), and $m_k$, $m_k'$ are 
some integers determined by $(L,\wB,D)$ 
(see Lemma \ref{lem:decat} and Definition \ref{def:monoidal mutation}
for more precise definition).

We say that $\shc$ is a {\it monoidal categorification} of a
skew-symmetric quantum cluster algebra $A$ if {\rm (i)}
the Grothendieck ring of $\shc$ is isomorphic to
$A$, {\rm (ii)} there exists a quantum monoidal seed
$\seed= (\{ M_i\}_{i \in \K}, L,\widetilde{B},D)$ in $\shc$ such that,
for some $m_i \in \frac{1}{2}\Z \ (i \in \K)$,
$[\seed] \seteq (\{ q^{m_i}[M_i]\}_{i \in \K}, L,\widetilde{B})$
is an initial quantum seed of $A$ {\rm (iii)} $\seed$ admits
successive mutations in all directions.

The existence of monoidal category $\shc$ which provides a
monoidal categorification of quantum cluster algebra $A$ implies the
followings: \eqn&&\parbox{80ex}{ \begin{enumerate}
\item[{\rm (QM1)}] Every quantum cluster monomial corresponds to the isomorphism class of a real simple object of $\shc$. In particular,
the set of quantum cluster monomials is $\Z[q^{\pm 1/2}]$-linearly independent.
\item[{\rm (QM2)}] The quantum positivity conjecture holds for $A$.
\end{enumerate}
} \eneqn

\medskip
In the last section of this
paper, we study the case when the category $\shc$ is a full subcategory
of $R \gmod$  which is stable under taking convolution products, subquotients,
extensions and grading shifts. Then the category has a natural
decomposition $\shc=\bigoplus_{\beta \in
\rtl^-}\shc_\beta$ for $\shc_\beta=\shc \cap
R(-\beta) \gmod$.

We say that a pair
$(\{ M_i \}_{i \in \K}, \wB)$ is {\em admissible} if {\rm (i)} $\{ M_i
\}_{i \in \K}$ is a family of real simple {\it self-dual} objects in
$\shc$ which commute with each other, {\rm (ii)} $\wB$ is an
integer-valued $\K\times \K_{{\rm ex}}$-matrix with skew-symmetric
part, {\rm (iii)} for any $k \in \K_{{\rm ex}}$, there exists a self-dual
simple object $M_k'$ in $\shc$ such that there is an exact
sequence in $\shc$
 \eqn
&&0 \to q \sodot_{b_{ik} >0} M_i^{\snconv b_{ik}} \to
q^{\tLa(M_k,M_k')} M_k \conv M_k' \to
 \sodot_{b_{ik} <0} M_i^{\snconv (-b_{ik})} \to 0,
 \eneqn
and $M_k'$ commutes with $M_i$  for any  $i \neq k$ (Definition \ref{def:admissible}).

The main theorem of our paper is the following (Theorem \ref{th:main},
Corollary \ref{cor:main}):

\begin{maintheorem}
Set $\La\seteq (\La(M_i,M_j))_{i,j\in \K}$, and let $\seed=(\{M_i\}_{i
\in \K},-\La,\wB,D)$ be a quantum monoidal seed in $\shc$.
 We assume further the following condition{\rm:}
\begin{itemize}
\item The Grothendieck ring of $\shc$ is isomorphic 
to the quantum cluster algebra $A$ associated to
the initial quantum seed $[\seed]\seteq \bl \{q^{-\frac{1}{4}(d_i,d_i)}
[M_i]\}_{i \in \K},-\La, \wB\br$ with $d_i:=\wt(M_i)$.
\end{itemize}
If the pair $(\{M_i\}_{i \in \K},\widetilde B)$ is admissible, then the
category $\shc$ is a monoidal categorification of the quantum
cluster algebra $A$. 
\end{maintheorem}

When the base field of the symmetric KLR algebra is of characteristic $0$,
Main Theorem, along with 
Theorem~\ref{thm:categorification 2} due to \cite{VV09, R11}, implies the
quantization conjecture: \eqn&&\parbox{80ex}{
\begin{enumerate}
\item[{\rm (QM3)}] The set of cluster monomials of $A$ is contained in the upper global basis $\B^\upper$.
\end{enumerate}
} \eneqn

In the forthcoming paper, 
we will prove the existence of admissible quantum monoidal seeds in 
$\shc_w$. Hence, Main Theorem shows that the category $\shc_w$
provides a monoidal categorification
of the quantum cluster algebra $A_q(\mathfrak{n}(w))$, and (QM3)
holds for $A_q(\mathfrak{n}(w))$.

\medskip

The paper is organized as follows. In Section 1, we briefly review
some of basic materials on quantum groups and KLR algebras.
In Section 2, we continue the study in
\cite{KKKO14} about convolution products,  heads and socles of
$R$-modules. In Section 3, we prove the first conjecture of Leclerc
in \cite{L03}. In Section 4, we recall the definition of quantum
cluster algebras. In Section 5, we give the definitions of a monoidal
seed, a quantum monoidal seed, a monoidal categorification of a cluster
algebra and a monoidal categorification of a quantum cluster algebra.
In Section 6, we prove Main Theorem by using the results of
previous sections.

\bigskip

\noindent
{\bf Acknowledgements.} The authors would like to express their gratitude to
Peter McNamara who informed us his result.
They would also like to express their gratitude to
Bernard Leclerc and Yoshiyuki Kimura for many fruitful discussions.
The last two authors gratefully acknowledge the hospitality of
Research Institute for Mathematical Sciences, Kyoto University
during their visits in 2014.

\section{Quantum groups and KLR algebras}
\subsection{Quantum groups and upper global bases} \label{subsec:qgroups}
Let $I$
be an index set. A \emph{Cartan datum} is a quintuple $(A,P,
\Pi,P^{\vee},\Pi^{\vee})$ consisting of
\begin{enumerate}[(a)]
\item an integer-valued matrix $A=(a_{ij})_{i,j \in I}$,
called the \emph{symmetrizable generalized Cartan matrix},
 which satisfies
\be[{\rm(1)}]
\item $a_{ii} = 2$ $(i \in I)$,
\item $a_{ij} \le 0 $ $(i \neq j)$,
\item there exists a diagonal matrix
$D=\text{diag} (\mathsf s_i \mid i \in I)$ such that $DA$ is
symmetric, and $\mathsf s_i$ are positive integers.
\end{enumerate}

\item a free abelian group $P$, called the \emph{weight lattice},
\item $\Pi= \{ \alpha_i \in P \mid \ i \in I \}$, called
the set of \emph{simple roots},
\item $P^{\vee}\seteq\Hom(P, \Z)$, called the \emph{co-weight lattice},
\item $\Pi^{\vee}= \{ h_i \ | \ i \in I \}\subset P^{\vee}$, called
the set of \emph{simple coroots},
satisfying the following properties:
\be[{\rm(1)}]
\item $\langle h_i,\alpha_j \rangle = a_{ij}$ for all $i,j \in I$,
\item $\Pi$ is linearly independent,
\item for each $i \in I$, there exists $\Lambda_i \in P$ such that
           $\langle h_j, \Lambda_i \rangle =\delta_{ij}$ for all $j \in I$.
\end{enumerate}
We call $\Lambda_i$ the \emph{fundamental weights}.
\end{enumerate}

\medskip
\noi
The free abelian group $\rootl\seteq\soplus_{i \in I} \Z \alpha_i$ is called the
\emph{root lattice}. Set $\rootl^{+}= \sum_{i \in I} \Z_{\ge 0}
\alpha_i\subset\rootl$ and $\rootl^{-}= \sum_{i \in I} \Z_{\le0}
\alpha_i\subset\rootl$. For $\beta=\sum_{i\in I}m_i\al_i\in\rootl$,
we set
$|\beta|=\sum_{i\in I}|m_i|$.

Set $\mathfrak{h}=\Q \otimes_\Z P^{\vee}$.
Then there exists a symmetric bilinear form $(\quad , \quad)$ on
$\mathfrak{h}^*$ satisfying
$$ (\alpha_i , \alpha_j) =\mathsf s_i a_{ij} \quad (i,j \in I)
\quad\text{and $\lan h_i,\lambda\ran=
\dfrac{2(\alpha_i,\lambda)}{(\alpha_i,\alpha_i)}$ for any $\lambda\in\mathfrak{h}^*$ and $i \in I$}.$$

Let $q$ be an indeterminate. For each $i \in I$, set $q_i = q^{\,\mathsf s_i}$.

\begin{definition} \label{def:qgroup}
The {\em quantum group} 
associated with a Cartan datum
$(A,P,\Pi,P^{\vee}, \Pi^{\vee})$ is the  algebra $\U$ over
$\mathbb Q(q)$ generated by $e_i,f_i$ $(i \in I)$ and
$q^{h}$ $(h \in P^\vee)$ satisfying the following relations:
\begin{equation*}
\begin{aligned}
& q^0=1,\ q^{h} q^{h'}=q^{h+h'} \ \ \text{for} \ h,h' \in P,\\
& q^{h}e_i q^{-h}= q^{\lan h, \alpha_i\ran} e_i, \ \
          \ q^{h}f_i q^{-h} = q^{-\lan h, \alpha_i\ran} f_i \ \ \text{for} \ h \in P^\vee, i \in
          I, \\
& e_if_j - f_je_i = \delta_{ij} \dfrac{K_i -K^{-1}_i}{q_i- q^{-1}_i
}, \ \ \mbox{ where } K_i=q^{\mathsf s_i h_i}, \\
& \sum^{1-a_{ij}}_{r=0} (-1)^r \left[\begin{matrix}1-a_{ij}
\\ r\\ \end{matrix} \right]_i e^{1-a_{ij}-r}_i
         e_j e^{r}_i =0 \quad \text{ if } i \ne j, \\
& \sum^{1-a_{ij}}_{r=0} (-1)^r \left[\begin{matrix}1-a_{ij}
\\ r\\ \end{matrix} \right]_i f^{1-a_{ij}-r}_if_j
        f^{r}_i=0 \quad \text{ if } i \ne j.
\end{aligned}
\end{equation*}
\end{definition}

Here, we set $[n]_i =\dfrac{ q^n_{i} - q^{-n}_{i} }{ q_{i} - q^{-1}_{i} },\quad
  [n]_i! = \prod^{n}_{k=1} [k]_i$ and
  $\left[\begin{matrix}m \\ n\\ \end{matrix} \right]_i= \dfrac{ [m]_i! }{[m-n]_i! [n]_i! }\;$
  for $i \in I$ and $m,n \in Z_{\ge 0}$ such that $m\ge n$.

Let $U_q^{+}(\g)$ (resp.\ $U_q^{-}(\g)$) be the subalgebra of
$U_q(\g)$ generated by $e_i$'s (resp.\ $f_i$'s), and let $U^0_q(\g)$
be the subalgebra of $U_q(\g)$ generated by $q^{h}$ $(h \in
P^{\vee})$. Then we have the \emph{triangular decomposition}
$$ U_q(\g) \simeq U^{-}_q(\g) \otimes U^{0}_q(\g) \otimes U^{+}_q(\g),$$
and the {\em weight space decomposition}
$$U_q(\g) = \bigoplus_{\beta \in \rootl} U_q(\g)_{\beta},$$
where $U_q(\g)_{\beta}\seteq\set{ x \in U_q(\g)}{\text{$q^{h}x q^{-h}
=q^{(h, \beta )}x$ for any $h \in P$}}$.

For any $i \in I$,
there exists a unique $\Q(q)$-linear endomorphism $e'_i$
of $U^-_q(\g)$  such that
\eqn
e'_i(f_j)=\delta_{i,j} \ (j \in I), &
 e'_i(xy) = (e'_i x) y + q_i^{\langle h_i, \beta \rangle} x (e'_iy) \ (x \in U^{-}_q(\g)_{\beta}, y \in U^-_q(\g)).
\eneqn

Recall that there exists a unique non-degenerate symmetric bilinear form $(\ ,\ )_{-}$ on $U^-_q(\g)$ such that
\eqn (\one,\one )=1, & (f_i u, v)_{-}=(u,e'_i v)_{-} \quad
\text{for $i \in I$, $u,v  \in U_q^-(\g)$}.
\eneqn

Let $\A= \Z[q, q^{-1}]$ and set
$$e_i^{(n)} = e_i^n / [n]_i!, \quad f_i^{(n)} =
f_i^n / [n]_i! \ \ (n \in \Z_{\ge 0}).$$
Let $U_{\A}^{-}(\g)$  be the
$\A$-subalgebra of $U^-_q(\g)$ generated by $f_i^{(n)}$  for $i\in I$, $n \in \Z_{\ge 0}$. It is called the {\em integral form} of $U^-_q(\g)$. The {\em dual integral form} $U^-_\A(\mathfrak g)^{\vee}$ of $U^-_q(\g)$ is
defined by $$U^-_\A(\mathfrak g)^{\vee} \seteq \{x \in U_q^-(\g) \,|\, (x, U_\A^-(\g))_{-}\subset \A\}.$$

Then $U^-_\A(\g)^\vee$ has an $\A$-algebra structure as a subalgebra of $U^-_q(\g)$.

For each $ i \in  I$, any element $x \in U_q^-(\g)$ can be written uniquely as
\eqn
x = \sum_{n \ge 0} f_i^{(n)} x_n\quad\text{ with $x_n \in \Ker(e'_i)$.}
\eneqn
We define the {\em Kashiwara operators} $\tilde e_i, \tilde f_i$ on
$U_q^-(\g)$ by
\eqn
\tilde e_i x = \sum_{n \ge 1} f_i^{(n-1)}x_n, & \tilde f_i x = \displaystyle\sum_{n \ge 0} f_i^{(n+1)}x_n,
\eneqn
and set
\eqn
&L(\infty) = \sum_{\ell \ge0, i_i, \ldots, i_\ell \in I}
\A_0 \tilde f_{i_1} \cdots \tilde f_{i_\ell} \cdot \one \subset U_q^-(\g),
\quad \overline{L(\infty)} = \set{\overline x}{ x \in L(\infty)}, \\
&B(\infty) = \set{\tilde f_{i_1} \cdots \tilde f_{i_\ell} \cdot \one \mod q L(\infty)}{\ell \ge0, i_i, \ldots, i_\ell \in I}  \subset L(\infty) / q L(\infty),
\eneqn
where $\A_0=\set{g \in \Q(q)}{g  \ \text{is regular at} \ q=0}$ and 
$- \col U_q(\g) \rightarrow U_q(\g)$ is the 
$\Q$-algebra involution  given by
 \eqn \overline{e_i} ={e_i}, \ \overline {f_i}=f_i, \ \overline{q^h}=q^{-h},  \ \text{and} \  \overline q =q^{-1}.
 \eneqn

The set $B(\infty)$ is a $\Q$-basis of $L(\infty) / q L(\infty)$ and the natural map
\eqn L(\infty) \cap \overline{L(\infty)} \cap U^-_\A(\g) \rightarrow L(\infty) / q L(\infty)\eneqn
is a $\Q$-linear isomorphism.
Let us denote the inverse of the above isomorphism by $G^{\low}$.
Then the set
$$\B^\low \seteq  \set{G^\low(b)}{b \in B(\infty)}$$
forms an $\A$-basis of $U_\A^-(\g)$ and is  called  the  {\em lower global basis} of $U_q^-(\g)$.

For each $b \in B(\infty)$, define $G^{\upper}(b) \in U^-_q(\g)$ as the element satisfying
$$(G^{\upper}(b), G^{\low}(b'))_{-} =\delta_{b,b'} $$ for $b' \in B(\infty)$.
The set
$$\B^\upper \seteq  \set{G^\upper(b)}{b \in B(\infty)}$$
forms an $\A$-basis of $U_\A^-(\g)^\vee$ and is called the {\em upper global basis} of $U_q^-(\g)$.

\subsection{\KLRs\ }
 \hfill

Now we recall the definition of Khovanov-Lauda-Rouquier algebra or quiver Hecke algebras
(hereafter, we abbreviate it 
as KLR algebras) associated with a given
Cartan datum $(A, P, \Pi, P^{\vee}, \Pi^{\vee})$.

Let $\cor$ be a base field.
For $i,j\in I$ such that $i\not=j$, set
$$S_{i,j}=\set{(p,q)\in\Z_{\ge0}^2}{(\al_i , \al_i)p+(\al_j , \al_j)q=-2(\al_i , \al_j)}.
$$
Let us take  a family of polynomials $(Q_{ij})_{i,j\in I}$ in $\cor[u,v]$
which are of the form
\eq
&&\parbox{68ex}{
$
Q_{ij}(u,v) = \begin{cases}\hs{5ex} 0 \ \ & \text{if $i=j$,} \\[1.5ex]
\sum\limits_{(p,q)\in S_{i,j}}
t_{i,j;p,q} u^p v^q\quad& \text{if $i \neq j$}
\end{cases}
$

\vs{1ex}
\hs{8ex}\parbox{65ex}{
with $t_{i,j;p,q}\in\cor$ such that
$Q_{i,j}(u,v)=Q_{j,i}(v,u)$ and
$t_{i,j:-a_{ij},0} \in \cor^{\times}$.}
}
\label{eq:Q}
\eneq

We denote by
$\sym_{n} = \langle s_1, \ldots, s_{n-1} \rangle$ the symmetric group
on $n$ letters, where $s_i\seteq (i, i+1)$ is the transposition of $i$ and $i+1$.
Then $\sym_n$ acts on $I^n$ by place permutations.

For $n \in \Z_{\ge 0}$ and $\beta \in \rootl^+$ such that $|\beta| = n$, we set
$$I^{\beta} = \set{\nu = (\nu_1, \ldots, \nu_n) \in I^{n}}%
{ \alpha_{\nu_1} + \cdots + \alpha_{\nu_n} = \beta }.$$

\begin{definition}
For $\beta \in \rootl^+$ with $|\beta|=n$, the {\em
Khovanov-Lauda-Rouquier algebra}  $R(\beta)$  at $\beta$ associated
with a Cartan datum $(A,P, \Pi,P^{\vee},\Pi^{\vee})$ and a matrix
$(Q_{ij})_{i,j \in I}$ is the  algebra over $\cor$
generated by the elements $\{ e(\nu) \}_{\nu \in  I^{\beta}}$, $
\{x_k \}_{1 \le k \le n}$, $\{ \tau_m \}_{1 \le m \le n-1}$
satisfying the following defining relations:
\allowdisplaybreaks[4]
\begin{equation*} \label{eq:KLR}
\begin{aligned}
& e(\nu) e(\nu') = \delta_{\nu, \nu'} e(\nu), \ \
\sum_{\nu \in  I^{\beta} } e(\nu) = 1, \\
& x_{k} x_{m} = x_{m} x_{k}, \ \ x_{k} e(\nu) = e(\nu) x_{k}, \\
& \tau_{m} e(\nu) = e(s_{m}(\nu)) \tau_{m}, \ \ \tau_{k} \tau_{m} =
\tau_{m} \tau_{k} \ \ \text{if} \ |k-m|>1, \\
& \tau_{k}^2 e(\nu) = Q_{\nu_{k}, \nu_{k+1}} (x_{k}, x_{k+1})
e(\nu), \\
& (\tau_{k} x_{m} - x_{s_k(m)} \tau_{k}) e(\nu) = \begin{cases}
-e(\nu) \ \ & \text{if} \ m=k, \nu_{k} = \nu_{k+1}, \\
e(\nu) \ \ & \text{if} \ m=k+1, \nu_{k}=\nu_{k+1}, \\
0 \ \ & \text{otherwise},
\end{cases} \\
& (\tau_{k+1} \tau_{k} \tau_{k+1}-\tau_{k} \tau_{k+1} \tau_{k}) e(\nu)\\
& =\begin{cases} \dfrac{Q_{\nu_{k}, \nu_{k+1}}(x_{k},
x_{k+1}) - Q_{\nu_{k}, \nu_{k+1}}(x_{k+2}, x_{k+1})} {x_{k} -
x_{k+2}}e(\nu) \ \ & \text{if} \
\nu_{k} = \nu_{k+2}, \\
0 \ \ & \text{otherwise}.
\end{cases}
\end{aligned}
\end{equation*}
\end{definition}

The above relations are homogeneous provided with
\begin{equation*} \label{eq:Z-grading}
\deg e(\nu) =0, \quad \deg\, x_{k} e(\nu) = (\alpha_{\nu_k}
, \alpha_{\nu_k}), \quad\deg\, \tau_{l} e(\nu) = -
(\alpha_{\nu_l} , \alpha_{\nu_{l+1}}),
\end{equation*}
and hence $R( \beta )$ is a $\Z$-graded algebra.

 For a graded $R(\beta)$-module $M=\bigoplus_{k \in \Z} M_k$, we define
$qM =\bigoplus_{k \in \Z} (qM)_k$, where
 \begin{align*}
 (qM)_k = M_{k-1} & \ (k \in \Z).
 \end{align*}
We call $q$ the \emph{grading shift functor} on the category of
graded $R(\beta)$-modules.

If $M$ is an $R(\beta)$-module, then we set $\wt(M)=-\beta \in \rootl^-$
and call it the {\em weight} of $M$.

\smallskip
We denote by $R(\beta) \Mod$ 
the category of $R(\beta)$-modules,
and by $R(\beta) \smod$
the full subcategory of $R(\beta)\Mod$ consisting of modules $M$
such that $M$ are finite-dimensional over $\cor$, 
and the actions of the $x_k$'s on $M$ are nilpotent.

Similarly, we denote by $R(\beta)\gMod$ and by $R(\beta)\gmod$
the category of graded $R(\beta)$-modules and
the category of graded $R(\beta)$-modules which are finite-dimensional over $\cor$, respectively.
We set 
$$
R\gmod=\soplus_{\beta\in\rtl^+}R(\beta)\gmod\quad\text{and}\quad
R\smod=\soplus_{\beta\in\rtl^+}R(\beta)\smod.
$$

\medskip
For $\beta, \gamma \in \rootl^+$ with $|\beta|=m$, $|\gamma|= n$,
 set
$$e(\beta,\gamma)=\displaystyle\sum_{\substack{\nu \in I^{\beta+\gamma}, \\ (\nu_1, \ldots ,\nu_m) \in I^{\beta}}} e(\nu) \in R(\beta+\gamma). $$
Then $e(\beta,\gamma)$ is an idempotent.
Let
\eq R( \beta)\tens R( \gamma  )\to e(\beta,\gamma)R( \beta+\gamma)e(\beta,\gamma) \label{eq:embedding}
\eneq
be the $\cor$-algebra homomorphism given by
$e(\mu)\tens e(\nu)\mapsto e(\mu*\nu)$ ($\mu\in I^{\beta}$
and $\nu\in I^{\gamma}$)
$x_k\tens 1\mapsto x_ke(\beta,\gamma)$ ($1\le k\le m$),
$1\tens x_k\mapsto x_{m+k}e(\beta,\gamma)$ ($1\le k\le n$),
$\tau_k\tens 1\mapsto \tau_ke(\beta,\gamma)$ ($1\le k<m$),
$1\tens \tau_k\mapsto \tau_{m+k}e(\beta,\gamma)$ ($1\le k<n$).
Here $\mu*\nu$ is the concatenation of $\mu$ and $\nu$;
i.e., $\mu*\nu=(\mu_1,\ldots,\mu_m,\nu_1,\ldots,\nu_n)$.

\medskip
For an $R(\beta)$-module $M$ and an $R(\gamma)$-module $N$,
we define the \emph{convolution product}
$M\conv N$ by
$$M\conv N=R(\beta + \gamma) e(\beta,\gamma)
\tens_{R(\beta )\otimes R( \gamma)}(M\otimes N). $$

For $M \in R( \beta) \smod$, the dual space
$$M^* \seteq \Hom_{\cor}(M, \cor)$$
admits an $R(\beta)$-module structure via
\begin{align*}
(r \cdot  f)(u) \seteq f(\psi(r) u) \quad (r \in R( \beta), \ u \in M),
\end{align*}
where $\psi$ denotes the $\cor$-algebra anti-involution on $R(\beta
)$ which fixes the generators $e(\nu)$, $x_m$ and $\tau_k$ for $\nu
\in I^{\beta}, 1 \leq m \leq  |\beta|$ and $1 \leq k<|\beta|$.

It is known that (see \cite[Theorem 2.2 (2)]{LV11})
\eq
&&(M_1 \conv M_2)^* \simeq q^{(\beta,\gamma)}
(M_2^* \conv M_1 ^*)\label{eq:dualconv}
\eneq
for any $M_1 \in R(\beta) \gmod$ and $M_2 \in R(\gamma) \gmod$.

A simple module $M$ in $R \gmod$ is called \emph{self-dual} if $M^* \simeq M$.
Every simple module is isomorphic to a grading shift of a self-dual simple module (\cite[\S.3.2]{KL09}).

\medskip
Let us denote by $K(R\gmod)$ the  Grothendieck group of $R\gmod$.
Then,
$K(R\gmod)$
is an algebra over  $\A\seteq\Z[q^{\pm1}]$ with
 the multiplication induced by the convolution product and
the $\A$-action induced by the  grading shift functor $q$.

 In \cite{KL09, R08}, it is shown that
\KLRs\ \emph{categorify} the negative half of the corresponding quantum group. More precisely, we have the following theorem.

\begin{theorem}[{\cite{KL09, R08}}] \label{Thm:categorification}
 For a given  Cartan datum $(A,P, \Pi,P^{\vee},\Pi^{\vee})$,
we take a parameter matrix $(Q_{ij})_{i,j \in \K}$ satisfying the conditions in \eqref{eq:Q}, and let $U_q(\mathfrak g)$ and $R(\beta) $ be
the associated quantum group   and the \KLRs, respectively.
Then there exists an $\A$-algebra isomorphism
\begin{align}
  U^-_\A(\mathfrak g)^{\vee} \simeq K(R\gmod).
\label{eq:KLRU}
\end{align}
\end{theorem}

 \KLRs\ also categorify the upper global bases.

\Def We say that the \KLR\ $R$ is {\em symmetric} if
$Q_{i,j}(u,v)$ is a polynomial in $u-v$ for all $i,j\in I$.
\edf
In particular, the generalized Cartan matrix $A$ is symmetric.

\begin{theorem} [\cite{VV09, R11}] \label{thm:categorification 2}
Assume that \KLR\ $R$ is symmetric and
the base field $\cor$ is of characteristic zero.
 Then under the isomorphism \eqref{eq:KLRU}
in {\rm Theorem \ref{Thm:categorification}},
the upper global basis  corresponds to
the set of the isomorphism classes of self-dual simple $R$-modules.
\end{theorem}

\subsection{R-matrices for Khovanov-Lauda-Rouquier algebras} \hfill

  For $|\beta|=n$ and $1\le a<n$,  we define $\vphi_a\in R( \beta)$
by
\eq&&\ba{l}
  \vphi_a e(\nu)=
\begin{cases}
  \bl\tau_ax_a-x_{a}\tau_a\br e(\nu)
  & \text{if $\nu_a=\nu_{a+1}$,} \\[2ex]
\tau_ae(\nu)& \text{otherwise.}
\end{cases}
\label{def:int} \ea \eneq
They are called the {\em intertwiners}. Since $\{\vphi_a\}_{1\le a<n}$
satisfies the braid relation, $\vphi_w\seteq\vphi_{i_1}\cdots\vphi_{i_\ell}$
does not depend on the choice of reduced expression $w=s_{i_1}a\cdots s_{i_\ell}$.

For $m,n\in\Z_{\ge0}$,
let us denote by $w[{m,n}]$ the element of $\sym_{m+n}$  defined by
\eq
&&w[{m,n}](k)=\begin{cases}k+n&\text{if $1\le k\le m$,}\\
k-m&\text{if $m<k\le m+n$.}\end{cases}
\eneq

Let $\beta,\gamma\in \rtl^+$ with $|\beta|=m$, $|\gamma|=n$,
 and let
$M$ be an $R(\beta)$-module and $N$ an $R(\gamma)$-module. Then the
map $M\tens N\to N\conv M$ given by
$u\tens v\longmapsto \vphi_{w[n,m]}(v\tens u)$
is  $R( \beta)\tens R(\gamma )$-linear, and  hence it extends
to an $R( \beta +\gamma)$-module  homomorphism
\eq &&R_{M,N}\col
M\conv N\To N\conv M. \eneq

\medskip
Assume that the \KLR\ $R(\beta)$ is symmetric.
Let $z$ be an indeterminate  which is homogeneous of degree $2$, and
let $\psi_z$ be the graded algebra homomorphism
\eqn
&&\psi_z\col R( \beta )\to \cor[z]\tens R( \beta )
\eneqn
given by
$$\psi_z(x_k)=x_k+z,\quad\psi_z(\tau_k)=\tau_k, \quad\psi_z(e(\nu))=e(\nu).$$

For an $R( \beta )$-module $M$, we denote by $M_z$
the $\bl\cor[z]\tens R( \beta )\br$-module
$\cor[z]\tens M$ with the action of $R( \beta )$ twisted by $\psi_z$.
Namely,
\eq&&
\ba{l}
e(\nu)(a\tens u)=a\tens e(\nu)u,\\[1ex]
x_k(a\tens u)=(za)\tens u+a\tens (x_ku),\\[1ex]
\tau_k(a\tens u)=a\tens(\tau_k u)
\ea
\label{eq:spt1}
\eneq
for  $\nu\in I^\beta$,  $a\in \cor[z]$ and $u\in M$.
 For $u\in M$, we sometimes denote by $u_z$ the corresponding element $1\tens u$ of
the $R( \beta )$-module  $M_z$.

For a non-zero $M\in R(\beta)\smod$
and a non-zero $N\in R(\gamma)\smod$,
\eq&&\parbox{70ex}{%
let $s$ be the order of zero of $R_{M_z,N}\col  M_z\conv
N\To N\conv M_z$;
i.e., the largest non-negative integer such that the image of
$R_{M_z,N}$ is contained in $z^s(N\conv M_z)$.}
\label{def:s} \eneq
Note that such an $s$ exists because
$R_{M_z,N}$ does not vanish (\cite[Proposition 1.4.4 (iii)]{K^3}).
We denote by $\Rm_{M_z,N}$
the morphism $z^{-s}R_{M_z,N}$.

\Def Assume that $R(\beta)$ is symmetric.
For a non-zero $M\in R(\beta)\smod$
and a non-zero $N\in R(\gamma)\smod$,
let $s$ be an integer as in \eqref{def:s}.
We define $$\rmat{M,N}\col M\conv N\to N\conv M$$
by  $$\rmat{M,N} = \Rm_{M_z,N}\vert_{z=0},$$
and call it the {\em renormalized R-matrix}.
\edf

By the definition, the renormalized R-matrix
$\rmat{M,N}$ never vanishes.

We define also
$$\rmat{N,M}\col N\conv M\to M\conv N$$
by  $$\rmat{N,M} = \bl (-z)^{-t}R_{N,M_z}\br\vert_{z=0},$$
where $t$ is the order of zero of
$R_{N,M_z}$.

If $R(\beta)$ and $R(\gamma)$ are symmetric, then
$s$ coincides with the multiplicity of zero of
$R_{M,N_z}$, and $\bl z^{-s}R_{M_z,N}\br\vert_{z=0}=\bl (-z)^{-s}R_{M,N_z}\br\vert_{z=0}$
(see, \cite[(1.11)]{KKKO14}).

By the construction, if the composition
$(N_1 \conv \rmat{M,N_2}) \circ (\rmat{M,N_1} \conv N_2)$
for $M,N_1,N_2 \in R \smod$
doesn't vanish then it is equal to $\rmat{M,N_1 \circ N_2}$.

\Def
A simple $R(\beta)$-module $M$ is called \emph{real} if $M \conv M$ is simple.
\edf

\Th[{\cite[Theorem 3.2]{KKKO14}}]
Let $\beta, \gamma \in \rootl^+$ and assume that $R(\beta)$ is symmetric.
Let $M$ be a real simple module in 
$R(\beta) \smod$ and
$N$ a simple module in $R(\gamma) \smod$.
Then
\bnum
  \item $M \conv N$ and $N \conv M$ have simple socles and simple heads.
  \item Moreover, $\Im(\rmat{M,N})$ is equal to the head of $M \conv N$
and socle of $N \conv M$.
Similarly,
$\Im(\rmat{N,M})$ is equal to the head of $N \conv M$
and socle of $M \conv N$.
\end{enumerate}
\enth

\medskip

\section{Simplicity of heads and socles of convolution products}

{\em In the rest of this paper, 
we assume that $R(\beta)$ is symmetric} for any $\beta\in\rtl^+$,
i.e., $Q_{ij}(u,v)$ is a function in $u-v$ for any $i,j\in I$.

We also work always in the category of graded modules.
For the sake of simplicity, {\em
we simply say that $M$ is an $R$-module instead of saying
that $M$ is a graded $R(\beta)$-module for $\beta\in\rtl^+$.}
We also sometimes ignore grading shifts if there is no afraid of confusion.
Hence, for $R$-modules $M$ and $N$, {\em we sometimes say that
$f\col M\to N$ is a homomorphism if
$f\col q^aM\to N$ is a morphism in $R\gmod$ for some $a\in\Z$.}
If we want to emphasize that $f\col q^aM\to N$ is a morphism in $R\gmod$,
we say so.

\Prop
Let $M$ be a real simple module, and let $N$ be a module with a simple socle.
If the following diagram
$$\xymatrix@C=12ex
{\wb{\soc(N)\conv M}\ar[r]^{\rmat{\soc(N),M}}\ar@{>->}[d]
&\wb{M\conv\soc(N)}\ar@{>->}[d]\\
N\conv M\ar[r]^{\rmat{N,M}}&M\conv N
}
$$ commutes up to a constant multiple, then
$\soc\bl M\conv\soc(N)\br$ is
equal to the socle of $M\conv N$.
In particular $M\conv N$ has a simple socle.
\enprop
\Proof
Let $S$ be an arbitrary simple submodule of $M\conv N$.
Then we have a commutative diagram
$$\xymatrix@C=12ex
{\wb{S\conv M}\ar[r]\ar@{>->}[d]
&\wb{M\conv S}\ar@{>->}[d]\\
M\conv N\conv M\ar[r]^{M\circ\rmat{N,M}}&M\conv M\conv N.
}$$
Hence $S\conv M\subset M\conv (\rmat{N,M})^{-1}(S)$.
Hence there exists a submodule $K$ of $N$ such that
$S\subset M\conv K$ and $K\conv M\subset (\rmat{N,M})^{-1}(S)$.
Hence $K\not=0$ and $\soc(N)\subset K$.
Hence $\rmat{N,M}\bl\soc(N)\conv M\br
\subset \rmat{N,M}\bl K\conv M\br\subset S$.
Since $\rmat{N,M}\bl\soc(N)\conv M\br$ is non-zero by the assumption, we have
$\rmat{N,M}\bl\soc(N)\conv M\br=S$.
Hence we obtain the desired result.
\QED

The following is a dual form of the preceding proposition.
\Prop\label{prop:simpleH}
Let $M$ be a real simple module.
Let $N$ be a module with a simple head.
If the following diagram
$$\xymatrix@C=12ex
{{M\conv N}\ar[r]^{\rmat{M,N}}\ar@{->>}[d]
&{N\conv M}\ar@{->>}[d]\\
M\conv \hd(N)\ar[r]^{\rmat{M,\hd(N)}}&\hd(N)\conv M
}
$$ commutes up to a constant multiple, then
$M\hconv \hd(N)$ is
equal to the simple head of $M\conv N$.
\enprop

\Def
For non-zero $M,N \in R \gmod$, we denote by $\Lambda(M,N)$ the homogeneous degree of the R-matrix
$\rmat{M,N}$. 
\edf

Hence 
\eqn
&&\Rm_{M_z,N}:M_z \conv N \to   q^{-\Lambda(M,N)} N \conv M_z \quad \text{and}\\
&&\rmat{M,N}:M\conv N \to q^{-\Lambda(M,N)} N \conv M
\eneqn
are morphisms in $R\gMod$ and in $R \gmod$, respectively.

\Def
For simple $R$-modules $M$ and $N$, we denote by $M \hconv N$ the head of
$M \conv N$.
\edf

\Lemma \label{lem:tLa even}
 For non-zero $R$-modules $M$ and $N$, we have
\eqn
\Lambda(M,N) \equiv \bl\wt(M),\wt(N)\br \mod 2.
\eneqn
\enlemma
\Proof
Set $\beta\seteq-\wt(M)$ and $\gamma\seteq-\wt(N)$.
By \cite[(1.3.3)]{K^3}, the homogeneous degree of  $R_{M_z,N}$ is $-(\beta,\gamma)+2\inp{\beta,\gamma}$, where
$\inp{{\scbul,\scbul}}$ is the symmetric bilinear form on $\rootl$ given by
$\inp{\al_i,\al_j}=\delta_{ij}$.
Hence  $\Rm_{M_z,N}=z^{-s}R_{M_z,N}$  has degree $-(\beta,\gamma)+2\inp{\beta,\gamma}-2s$.
\QED

\Def
For non-zero $R$-modules $M$ and $N$, we set
\eqn
\tL(M,N) \seteq \frac{1}{2}\bl\Lambda(M,N) +(\wt(M), \wt(N))\br.
\eneqn
\edf

\Lemma\label{lem:selfdual}
Let $M$ and $N$ be self-dual simple modules.
If one of them is real, then
\eqn
q^{\tL(M,N)} M \hconv N
\eneqn
is a self-dual simple module.
\enlemma
\Proof
Set $\beta=\wt(M)$ and $\gamma=\wt(N)$.
Set $M \hconv N = q^c L$ for some self-dual simple module $L$ and some $c \in \Z$.
Then we have
\eqn
M \conv N \epito q^c L \monoto q^{-\Lambda(M,N)} N \conv M.
\eneqn
Taking dual, we obtain
\eqn
q^{\Lambda(M,N)+(\beta,\gamma)} M \conv N \epito q^{-c} L \monoto q^{(\beta,\gamma)} N \conv M.
\eneqn
Hence we have $c=-c-\Lambda(M,N)-(\beta,\gamma)$, which implies $c= -\tL(M,N)$.
\QED

\Lemma
Let $M$ be a simple module and let $N_1, N_2$ be non-zero modules.
Then the composition
\eqn
M \conv N_1 \conv N_2 \To[\rmat{M,N_1} \circ N_2]
N_1  \conv M \conv N_2 \To[N_1 \circ\; \rmat{M,N_2}]
N_1 \conv N_2 \conv M
\eneqn
coincides with $\rmat{M,N_1 \circ N_2}$,
 and the composition
\eqn
N_1 \conv N_2 \conv M \To[N_1\circ\; \rmat{N_2,M}]
N_1  \conv M \conv N_2 \To[\rmat{N_1,M} \circ N_2]
 M \conv N_1 \conv N_2
\eneqn
coincides with $\rmat{N_1 \circ N_2,N}$.

In particular, we have

\eqn
\Lambda(M, N_1\conv N_2) = \Lambda(M,N_1)+\Lambda(M,N_2)
\eneqn
and
\eqn
\Lambda(N_1\conv N_2,M) = \Lambda(N_1,M)+\Lambda(N_2,M).
\eneqn
\enlemma
\Proof
It is enough to show that the compositions
$(N_1 \conv \rmat{M,N_2}) \circ (\rmat{M,N_1} \conv N_2)$ and
$ (\rmat{N_1,M} \conv N_2) \circ (N_1 \conv \rmat{N_2,M})$
 do not vanish.

Assume that $(N_1 \conv \rmat{M,N_2}) \circ (\rmat{M,N_1} \conv N_2)$ vanishes.
Then we have
\eqn
\Im(\rmat{M,N_1}) \circ N_2 \subset N_1 \conv \Ker(\rmat{M,N_2}).
\eneqn
By \cite[Lemma 3.1]{KKKO14}, there is a submodule $Z$ of $M$ such that
\eqn
\Im(\rmat{M,N_1}) \subset N_1 \conv Z \quad \text{and} \ Z \conv N_2 \subset \Ker(\rmat{M,N_2}).
\eneqn
It contradicts  the assumption that $M$ is simple.

Similarly, one can show that $ (\rmat{N_1,M} \conv N_2) \circ (N_1 \conv \rmat{N_2,M})$ does not vanish.
\QED

\Lemma
Let $M$ and $N$ be simple $R$-modules.
Then we have
\bnum
   \item $\Lambda(M,N)+\Lambda(N,M) \in 2 \Z_{\ge 0}.$
   \item If $\Lambda(M,N)+\Lambda(N,M) =2m$ for some $m \in \Z_{\ge 0}$, then
  \eqn
\Rm_{M_z,N} \circ \Rm_{N,M_z}=z^m \id_{N \conv M_z} \quad
\text{and} \quad
\Rm_{N,M_z} \circ \Rm_{M_z,N}=z^m \id_{M_z \conv N}
  \eneqn
  up to constant multiples.
\end{enumerate}
\enlemma
\Proof
By \cite[Proposition 1.6.2]{K^3}, the morphism
\eqn
\Rm_{N,M_z} \circ \Rm_{M_z,N} : M_z \circ N \to 
M_z \circ N
\eneqn
is equal to $f(z) \id_{M_z \circ N}$ for some $0 \neq f(z) \in \cor[z]$.
Since $\Rm_{N,M_z} \circ \Rm_{M_z,N}$ is  homogeneous of degree $\Lambda(M,N)+\Lambda(N,M)$,
we have $f(z)=c z^{\frac{1}{2}(\Lambda(M,N)+\Lambda(N,M))}$ for some $c \in \cor^\times$.
\QED

\Def
For non-zero modules $M$ and $N$, we set
\eqn
\de(M,N) = \frac{1}{2}\bl\Lambda(M,N)+\Lambda(N,M)\br.
\eneqn
\edf
Note that if $M$ and $N$ are simple modules, then  we have
$\de(M,N) \in \Z_{\ge 0}$.
Note also that
if $M,N_1,N_2$ are simple modules, then we have
$\de(M,N_1 \conv N_2) =\de(M,N_1)+\de(M,N_2)$.
\Lemma[\cite{KKKO14}]
Let $M,N$ be simple modules and assume that one of them is real. Then the following conditions are equivalent.
\bna
  \item $\de(M,N)=0$.
  \item $\rmat{M,N}$ and $\rmat{N,M}$ are inverse to each other up to a constant multiple.
  \item $M \conv N$ and $N \conv M$ are isomorphic up to a grading shift.
  \item $M \hconv N$ and $N \hconv M$ are isomorphic up to a grading shift.
  \item $M \conv N$ is simple.
\end{enumerate}
\enlemma

\Def
Let $M,N$ be simple modules.
\bnum
  \item We say that $M$ and $N$ \emph{commute} if $\de(M,N)=0$.
  \item We say that $M$ and $N$ \emph{simply-linked} if $\de(M,N)=1$.
\end{enumerate}
\edf

\Prop  \label{prop:real simple}
Let $M_1,\ldots, M_r$ be a commuting family of real simple modules.
Then the convolution product
\eqn
M_{1} \conv \cdots \conv M_{r}
\eneqn
is a real simple module.
\enprop
\Proof We shall first show the simplicity of the convolutions.
By induction on $r$, we may assume that $ M_{2} \conv \cdots \conv M_{r}$ 
is a simple module.
Then we have
\eqn
\de(M_{1}, M_{2} \conv \cdots \conv M_{r})
=\sum_{s=2}^r \de(M_{1},M_{s})=0
\eneqn
so that
$M_{1} \conv \cdots \conv M_{r}$ is simple.

Since $(M_{1} \conv \cdots \conv M_{r}) \conv (M_{1} \conv \cdots \conv M_{r})$ is also simple,
$M_{1} \conv \cdots \conv M_{r}$ is real.
\QED

\Def Let $M_1,\ldots, M_m$ be real simple modules.
Assume that they are commuting with each other.
We set
$$\nconv_{1\le k\le m} M_k
\seteq q^{\sum_{1\le i<j\le m}\tL(M_i,M_j)}M_1\conv\cdots\conv M_m.$$
\edf

\Lemma\label{lem:multipro}
Let $M_1,\ldots, M_m$ be real simple modules commuting with each other.
Then for any $\sigma\in\sym_m$, we have
$$\nconv_{1\le k\le m} M_k\simeq\nconv_{1\le k\le m} M_{\sigma(k)}\quad\text{in $R\gmod$.}$$
Moreover, if the $M_k$'s are self-dual,
then so is
$\nconv_{1\le k\le m} M_k$.
\enlemma
\Proof
It follows from Lemma~\ref{lem:selfdual} and
$q^{\tL(M_i,M_j)}M_i\conv M_j\simeq
q^{\tL(M_j,M_i)}M_j\conv M_i$.
\QED

\Prop
Let $f:N_1 \to N_2$ be a morphism between non-zero modules 
$R$-modules $N_1, N_2 $
and let $M$ be a non-zero $R$-module.
\bnum
\item
If $\Lambda(M,N_1)=\Lambda(M,N_2)$, then the following diagram is commutative:
  \eqn
  \xymatrix@C=5em{
  M\conv N_1 \ar[r]^{\rmat{M,N_1}} \ar[d]_{M \circ f} &%
N_1 \conv M \ar[d]^{f \circ M} \\
  M\conv N_2 \ar[r]^{\rmat{M,N_2}}& %
N_2 \conv M.
  }\eneqn
\item If $\Lambda(M,N_1) < \Lambda(M,N_2)$, then the composition
  \eqn
  M \conv N_1 \To[M \circ f] M \conv N_2 \To[\rmat{M,N_2}] N_2 \conv M
  \eneqn
  vanishes.
 \item If $\Lambda(M,N_1) > \Lambda(M,N_2)$, then the composition
  \eqn
M \conv N_1 \To[\rmat{M,N_1}] N_1 \conv M \To[f \circ M] N_2 \conv M
  \eneqn
  vanishes.
 \item If $f$ is surjective, then we have
 \eqn
 \Lambda(M,N_1) \ge \Lambda(M,N_2) \quad \text{and} \quad
 \Lambda(N_1,M) \ge \Lambda(N_2,M)
 \eneqn
If $f$ is injective, then we have
\eqn
\Lambda(M,N_1) \le \Lambda(M,N_2) \quad \text{and}\quad
\Lambda(N_1,M) \le \Lambda(N_2,M)
\eneqn
\end{enumerate}

\enprop
\Proof

Let $s_i$  be the order of zero of $R_{M_z,N_i}$ for $i=1,2$.
Then we have  $\Lambda(M,N_1) - \Lambda(M,N_2)=2(s_2-s_1)$.

Set $m\seteq\min\{s_1,s_2\}$.
Then the following diagram is commutative:
  \eqn
  \xymatrix@C=10em{
  M_z\conv N_1 \ar[r]^{z^{-m}R_{M_z,N_1}} \ar[d]_{M_z \conv f} &  N_1 \conv M_z \ar[d]^{f \conv M_z} \\
  M_z\conv N_2 \ar[r]^{z^{-m}R_{M_z,N_2}}&  N_2 \conv M_z.
  }\eneqn

(i) If $s_1=s_2$, then by specializing $z=0$ in the above diagram,
we obtain the commutativity of the diagram in (i).

(ii) If $s_1 > s_2$, then we have
\eqn
z^{-m}R_{M_z,N_1}=z^{s_1-m} \bl z^{-s_1}R_{M_z,N_1}\br
\eneqn
so that $z^{-m}R_{M_z,N_1}|_{z=0}$ vanishes.
Hence
we have
\eqn
\rmat{M,N_2} \circ (M\conv f) = z^{-m}R_{M_z,N_2}|_{z=0} \circ (M\conv f) =0,
\eneqn
as desired. In particular, $f$ is not surjective.

(iii) Similarly, if $s_1 < s_2$, then we have
$
(f \conv M) \circ \rmat{M,N_1} =0,
$
and $f$ is not injective.

(iv) The statements for $\Lambda(M,N_1)$ and  $\Lambda(N_1, M)$ follow  
from (ii) and (iii).
The other statements can be shown in a similar way.
\QED

\Prop
Let $L$, $M$ and $N$ be simple modules.
Then we have
\eq&&\ba{l}
\La(L,S)\le\La(L,M)+\La(L,N),\ 
\La(S,L)\le\La(M,L)+\La(N,L)\ \text{and}\\[1ex]
\de(S,L)\le\de(M,L)+\de(N,L)\ea\label{eq:LMN<}
\eneq
for any subquotient $S$ of $M\conv N$.
Moreover, when $L$ is real, the following conditions are equivalent.
\bna
\item $L$ commutes with $M$ and $N$.
\item Any simple subquotient $S$ of $M\conv N$ commutes with 
$L$ and satisfies $\La(L,S)=\La(L,M)+\La(L,N)$.
\item Any simple subquotient $S$ of $M\conv N$ commutes with $L$ and satisfies
$\La(S,L)=\La(M,L)+\La(N,L)$.
\ee
\enprop
\Proof
The inequalities \eqref{eq:LMN<} are consequences of the preceding proposition.
Let us show the equivalence
of (a)--(c).

Let $M\conv N=K_0\supset K_1\supset \cdots\supset  K_\ell\supset K_{\ell+1}=0$ be a Jordan-H\"older series of $M\conv N$.
Then the renormalized R-matrix $\Rm_{L_z,M\circ N}=(M\conv\Rm_{L_z,N})\circ(\Rm_{L_z,M}\conv N)
\cl L_z\conv M\conv N\to M\conv N\conv L_z$ 
is homogeneous of degree $\La(L,M)+\La(L,N)$ and it sends $L_z\conv K_k$ to
$K_k\conv L_z$ for any $k\in\Z$. Hence
$f\seteq\rmat{L,M\circ N}=\Rm_{L_z,M\circ N}\vert_{z=0}$ sends
$L\conv K_k$ to $K_k\conv L$.

\medskip
\noi
First assume (a). 
Then $f$ is an isomorphism.
Hence $f\vert_{L\circ K_k}\cl L\conv K_k\to K_k\conv L$ is injective.
By comparing their dimension, $f\vert_{L\circ K_k}$ is an isomorphism,
Hence $f\vert_{L\circ (K_k/K_{k+1})}$ is an isomorphism of homogeneous degree
$\La(L,M)+\La(L,N)$. 
Hence we obtain (b).

\medskip
\noi
Conversely, assume (b). 
Then, 
$\Rm_{L_z,M\circ N}\vert_{L_z\conv (K_k/K_{k+1})}$ and $\Rm_{L_z,K_k/K_{k+1}}$ 
have the same homogeneous degree, and hence they should coincide.
It implies that $f\vert_{L\circ(K_k/K_{k+1})}=\rmat{L,K_k/K_{k+1}}$ is an isomorphism 
for any $k$.
Therefore $f=(M\conv\rmat{L,N})\circ(\rmat{L,M}\conv N)$ is an isomorphism,
which implies that $\rmat{L,N}$ and $\rmat{L,M}$ are isomorphisms.
Thus we obtain (a).

\medskip
\noi
Similarly (a) and (c) are equivalent.
\QED

\Lemma \label{lem:LMN1}
Let $L$, $M$ and  $N$ be simple modules.
 We assume that $L$ is real and commutes with $M$.
Then the diagram
$$\xymatrix@C=10ex
{L\conv (M\conv N)\ar[r]^{\rmat{L,M\circ N}}\ar[d]
&(M\conv N)\conv L\ar[d]\\
L\conv (M\hconv N)\ar[r]^{\rmat{L,M\shconv N}}
&(M\hconv N)\conv L}
$$
commutes.
\enlemma
\Proof
Otherwise
the composition

$$
\xymatrix@C=8ex{
L\conv M\conv N\ar[r]^-\sim_-{\rmat{L,M}\circ N}
&M\conv L\conv N\ar[r]_-{M\circ\rmat{L,N}}&M\conv N\conv L
\ar[r]&
(M\hconv N)\conv L}
$$
vanishes.
Hence we have
$$M\circ \Im(\rmat{L,N})\subset \Ker(M\conv N\to M\hconv N)\conv L.$$
Hence there exists a submodule $K$ of $N$ such that
$$\text{$\Im(\rmat{L,N})\subset K\conv L$ and
$M\conv K\subset\Ker(M\conv N\to M\hconv N)$.}$$
The first inclusion implies $K\not=0$ and the second implies $K\not=N$,
which contradicts the simplicity of $N$.
\QED

The following lemma can be proved similarly.
\Lemma\label{lem:LMN2}
Let $L$, $M$ and  $N$ be simple modules.
We assume that $L$ is real and commutes with $N$.
Then the diagram
$$\xymatrix@C=10ex
{(M\conv N)\conv L\ar[r]^{\rmat{M\circ N,L}}\ar[d]
&L\conv(M\conv N)\ar[d]\\
(M\hconv N)\conv L\ar[r]^{\rmat{M\shconv N,L}}
&L\conv(M\hconv N)}
$$
 commutes. %
\enlemma

The following proposition follows from Lemma~\ref{lem:LMN1} and Lemma~\ref{lem:LMN2}.
\Prop
Let $L$, $M$ and  $N$ be simple modules. Assume that $L$ is real.
Then we have
\bnum
\item
If $L$ and $M$ commute, then
$$\La(L,M\hconv N)=\La(L,M)+\La(L,N).$$

\item
If $L$ and $N$ commute, then
$$\La(M\hconv N,L)=\La(M,L)+\La(N,L).$$
\ee
\enprop

\Prop \label{prop:real}
Let $X,Y,M$ and $N$ be simple $R$-modules.
Assume that there is  an exact sequence
\eqn
0 \to X \to M \conv N \to Y \to 0,
\eneqn
$X \conv N$ and $Y \conv N$ are simple and $X \conv N \not\simeq Y \conv N$ 
as ungraded modules.
Then $N$ is a real simple module.
\enprop
\Proof
Assume that $N$ is not real.
Then $N \conv N$ is reducible and we have $\rmat{N,N} \neq c \id_{N \conv N}$ for any $c \in \cor$.
Note that $N \conv N$
is of length $2$, because $M\conv N \conv N$ is of length $2$.

Let $S$ be a simple  submodule of $N \conv N$.
Consider an exact sequence
$$0\To X\conv N\To M\conv N\conv N\To Y\conv N\To0.$$
Then we have 
\eq (X \conv N) \cap (M\conv S) =0.
\label{eq:SN}
\eneq
Indeed, if $(X\conv N) \subset (M \conv S)$, then there exists a submodule $Z$ of $N$ such that
$X \subset M \conv Z$ and $Z \conv N \subset S$.
It contradicts the simplicity of $N$.
Thus \eqref{eq:SN} holds.

Note that \eqref{eq:SN} implies that
$$M\conv S\simeq Y\conv N.$$

\medskip\noi
(a)\ Assume first
that $N \conv N$ is semisimple so that $N \conv N=S \oplus S'$ for some simple submodule
 $S'$ of $N \conv N$.
 Then
$M\conv S\simeq Y\conv N\simeq M\conv S'$.
Hence $M\conv S\simeq X\conv N\simeq M\conv S'$.
Therefore we obtain $X \conv N \simeq Y \conv N,$ which is a contradiction.

\medskip\noi
(b)\   Assume that $N \conv N$ is not semisimple so that $S$ is a unique non-zero proper submodule of
  $N \conv N$ and $(N \conv N) / S$ is a unique non-zero proper quotient of $N \conv N$.
  Without loss of generality, we may assume that $\cor$ is algebraically closed.
  Let $x \in \cor$ be an eigenvalue of $\rmat{N,N}$.
 Since $\rmat{N,N} \not\in\cor\id_{N \conv N}$, we have
 $0  \subsetneq \Im(\rmat{N,N}-x \id_{N\conv N}) \subsetneq N\conv N$.
It follows that
\eqn
S = \Im(\rmat{N,N}-x \id_{N\conv N}) \simeq (N \conv N) /S,
\eneqn
and  hence we have an exact sequence
\eqn
0 \To M\conv S \To M \conv N \conv N \To M \conv \bl(N \conv N) /S\br  \To 0.
\eneqn
Since $M \conv N \conv N $ is of length $2$, we have
\eqn
X \conv N \simeq M \conv S \simeq M\conv\bl(N \conv N) /S\br\simeq  Y \conv N,
\eneqn
 which is a contradiction.
 \QED
 
 \Cor \label{cor:real}
Let $X,Y,N$ be simple $R$-modules and let $M$ be a real simple  $R$-module.
 If we have an exact sequence
\eqn
0 \to X \to M \conv N \to Y \to 0
\eneqn
and if $X \conv N$ and $Y \conv N$ are simple,
then $N$ is a real simple module.
\encor
\Proof
Since $M$ is real and $M \conv N$ is not simple, $X$ is not isomorphic to $Y$
as an ungraded module. 
It follows that $X \conv N$ is not isomorphic to $Y \conv N$, 
because $K(R\gmod)$ is a domain so that 
$[X\conv N]=q^m [Y \conv N]$ for some $m \in \Z$ implies $[X]=q^m [Y]$. 
Now the assertion follows from the above proposition.
\QED
The following lemmas are used later.
\Lemma \label{lem:convpower}
Let $M$ and $N$ be real simple modules.
Assume that $M\hconv N$ is real and commutes with $N$.
Then for any $n\in\Z_{\ge0}$, we have
$$M^{\circ n}\hconv N^{\circ n}\simeq (M\hconv N)^{\circ n}
\quad\text{as ungraded modules.}$$
\enlemma
\Proof
Set $L=M\hconv N$. Since $L^{\circ n}$ is simple, it is enough to give 
an epimorphism $M^{\circ n}\conv N^{\circ n}\epito L^{\circ n}$.
We shall show it by induction on $n$.
For $n>0$, we have
\eqn
&&M^{\circ n}\conv N^{\circ n}\simeq
M^{\circ(n-1)}\conv M\conv N\conv N^{\circ(n-1)}\\
&&\hs{10ex}\epito M^{\circ(n-1)}\conv L\conv  N^{\circ(n-1)}
\simeq M^{\circ(n-1)}\conv  N^{\circ(n-1)}\conv L
\epito L^{\circ(n-1)}\conv L.
\eneqn
\QED

\Cor\label{cor:crde}
Let $M$ and $N$ be simple module.
Assume that one of them is real.
If there is an exact sequence
$$0\to q^mX\To M\circ N\To q^n Y\To 0$$
for self-dual simple modules $X$, $Y$ and 
integers $m$, $n$, then we have
$$\de(M,N)=m-n\ge0.$$
\encor
\Proof
We may assume that $M$ and $N$ are self-dual without loss of generality.
Then we have
$n=-\tLa(M,N)$.
Since $q^mX\simeq q^{\La(N,M)}N\hconv M
\simeq q^{\La(N,M)-\tL(N,M)}\bl q^{\tL(N,M)}N\hconv M\br$, we have
$m=\La(N,M)-\tL(N,M)$.
Thus we obtain
$$m-n=\La(N,M)-\tL(N,M)+\tLa(M,N)=\de(M,N).$$
\QED

\section{Leclerc conjecture}
Recall that 
$R$ is assumed to be a symmetric KLR algebra
over a base field $\cor$.

\subsection{Leclerc conjecture}
\Th\label{th:leclerc}
Let $M$ and $N$ be simple modules.
We assume that $M$ is real.
Then we have
the equalities in the Grothendieck group $K(R\gmod)${\rm:}
\bnum
\item
$[M\conv N]=[M\hconv N]+\sum_{k}[S_k]$\\[.5ex]
with simple modules $S_k$ such that $\La(M,S_k)<\La(M,M\hconv N)=\La(M,N)$,
\item
$[M\conv N]=[q^{\La(N,M)}N\hconv M]+\sum_{k}[S_k]$\\[.5ex]
with simple modules $S_k$ such that $\La(S_k,M)<\La(N\hconv M,M)=\La(N,M)$,
\item
$[N\conv M]=[N\hconv M]+\sum_{k}[S_k]$\\[.5ex]
with simple modules $S_k$ such that $\La(S_k,M)<\La(N\hconv M,M)=\La(N,M)$, and
\item
$[N\conv M]=[q^{\La(M,N)}M\hconv N]+\sum_{k}[S_k]$\\[.5ex]
with simple modules $S_k$ such that $\La(M,S_k)<\La(M,M\hconv N)=\La(M,N)$.
\ee
In particular, $M\hconv N$ as well as $N\hconv M$ appears only once
in the Jordan-H\"older series of $M\conv N$ in $R\smod$.
\enth

The following result is an immediate consequence of this theorem.
\Cor\label{cor:compest}
Let $M$ and $N$ be simple modules.
We assume that one of them is real.
Assume that $M$ and $N$ do not commute,
Then we have
the equality in the Grothendieck group $K(R\gmod)$
$$[M\conv N]=[M\hconv N]+[q^{\La(N,M)}N\hconv M]+\sum_{k}[S_k]$$
with simple modules $S_k$.
Moreover we have
\bnum
\item
If $M$ is real, then we have
$\La(M, N\hconv M)<\La(M,N)$, $\La(M\hconv N,M)<\La(N,M)$ and
$\La(M, S_k)<\La(M,N)$, $\La(S_k,M)<\La(N,M)$.
\item
If $N$ is real, then we have
$\La(N, M\hconv N)<\La(N,M)$, $\La(N,M\hconv N)<\La(N,M)$ and
$\La(N, S_k)<\La(N,M)$, $\La(S_k,N)<\La(M,N)$.
\ee
\encor

\medskip

\Proof[Proof of Theorem~\ref{th:leclerc}]
We shall prove only (i).
The other statements are proved similarly.

Let us take a Jordan-H\"older series of $M\conv N$
$$M\conv N=K_0\supset K_1\supset \cdots\supset K_\ell\supset K_{\ell+1}=0.$$
Then we have
$K_0/K_1\simeq M\hconv N$.
Let us consider the
renormalized R-matrix
$\Rm_{M_z,M\circ N}=(M\conv\Rm_{M_z,N})\circ(\Rm_{M_z,M}\conv N)$
$$\xymatrix@C=10ex{
M_z\conv M\conv N\ar[r]^{\Rm_{M_z,M}\circ N}
&M\conv M_z\conv N\ar[r]^{M\circ\Rm_{M_z,N}}&M\conv N\conv M_z.}
$$
Then $\Rm_{M_z,M\circ N}$ sends $M_z\conv K_k$ to $K_k\conv M_z$ for any $k$.
Hence evaluating the above diagram at $z=0$, we obtain
$$\xymatrix@C=10ex{
M\conv M\conv N\ar[r]^{M\circ\rmat{M,N}}
&M\conv N\conv M\\
\wb{M\conv K_1}\ar[r]\ar@{^{(}->}[u] &\wb{K_1\conv M\,.}\ar@{^{(}->}[u]}
$$
Since $\Im(\rmat{M,N}\cl M\conv N\to N\conv M)\simeq (M\conv N)/K_1$,
we have
$\rmat{M,N}(K_1)=0$.
Hence, $\Rm_{M_z,M\circ N}$ sends $M_z\conv K_1$ to
$(K_1\conv M_z)\cap z\bl (M\conv N)\conv M_z\br
=z(K_1\conv M_z)$.
Thus $z^{-1}\Rm_{M_z,M\circ N}\vert_{M_z\circ K_1}$ is well defined.
Then it sends
$M_z\conv K_k$ to $K_k\conv M_z$ for $k\ge1$.
Thus we obtain an R-matrix
$$z^{-1}\Rm_{M_z, M\circ N}\vert_{M_z\circ(K_k/K_{k+1})}\cl
M_z\conv (K_k/K_{k+1})\to (K_k/K_{k+1})\conv M_z
\quad\text{for $1\le k\le \ell$.}$$
Hence we have
$$\Rm_{M_z, K_k/K_{k+1}}=z^{-s_k}z^{-1}\Rm_{M_z, M\circ N}\vert_{M_z\circ(K_k/K_{k+1})}$$
for some $s_k\in \Z_{\ge0}$.
Since the homogeneous degree of $\Rm_{M_z,M\circ N}$
is $\La(M,M\circ N)=\La(M,N)$, we obtain
$$\La(M,K_k/K_{k+1})=\La(M,N)-2(1+s_k)<\La(M,N).$$
\QED

The following theorem is an application of the above theorem.

\Th \label{thm:divisible}
Let $\phi$ be an element of the Grothendieck group $K(R\gmod)$
given by
\eqn
\phi = \sum_{b \in B(\infty)} a_b [L_b],
\eneqn
where $L_b$ is the self-dual simple module corresponding to $b$ and $a_b \in \Z[q^{\pm1}]$.
Let $A$ be a real simple module in $R\gmod$.
Assume that we have an equality 
\eqn
\phi [A] = q^l [A] \phi
\eneqn
in $K(R\gmod)$ for some $l \in \Z$.
Then $A$ commutes with $L_b$ and
\eqn
l= \La(A,L_b)
\eneqn
for  every $b \in B(\infty)$ such that $a_b \neq 0$.

\enth
\Proof
Note that we have
\eqn
\phi [A] = \sum_b a_b [L_b \conv A]
= \sum_b a_b ([L_b \hconv A] + \sum_k [S_{b,k}]) \quad \text{and} \\
q^l [A] \phi = q^l \sum_b a_b [A \conv L_b ]
= q^l \sum_b a_b (q^{\La(L_b,A)}[L_b \hconv A] + \sum_k [S^{b,k}]),
\eneqn
for some simple modules $S_{b,k}$ and $S^{b,k}$ satisfying
\eqn
\La(S_{b,k},A) < \La(L_b,A) \quad \text{and} \quad\La(S^{b,k},A) < \La(L_b,A).
\eneqn

We may assume that $\set{ b\in B(\infty)}{a_b \neq 0}\not=\emptyset$. Set
\eqn
t\seteq\displaystyle\max \set{ \La(L_b,A)}{ a_b \neq 0}.
\eneqn
By taking the classes of self-dual simple modules $S$ with $\La(S,A)=t$
in the expansions of $\phi [A]$ and $q^l [A] \phi$, we obtain
\eqn
\sum_{\La(L_b,A)=t} a_b [L_b \hconv A ]
=\sum_{\La(L_b,A)=t} q^l a_b q^{\La(L_b,A)} [L_b \hconv A ].
\eneqn
In particular, we have
$t=-l$.

Set
\eqn
t'\seteq\displaystyle\max \set{ \La(A,L_b)}{a_b \neq 0}.
\eneqn

Then, by a similar argument we have
$t'=l$.

It follows that

\eqn
0= t+t' \ge  \La(L_b,A) + \La(A,L_b)\ge0
\eneqn
for every $b$ such that $a_b \neq 0$. Hence $A$ and $L_b$ commute.

Since
\eqn \sum a_b ~ q^{\La(A,L_b)} [A\conv L_b]
 =\sum a_b [L_b \conv A]
=\phi [A]=q^l [A] \phi
 =q^l \sum a_b [A\conv L_b],
\eneqn
we have
\eqn
l=\La(A,L_b)
\eneqn
for any $b$ such that $a_b \neq 0$, as desired.
\QED

\Cor
Let $M$ and $N$ be simple modules. Assume that one of them is real.
If $[M]$ and $[N]$ q-commute \ro i.e., $[M][N]=q^n[N][M]$ for some $n\in\Z$\rf, then
$M$ and $N$ commute. In particular, $M\conv N$ is simple.
\encor
The following corollary is an immediate consequence of the corollary above
and Theorem~\ref{thm:categorification 2}.
\Cor
Assume that the generalized Cartan matrix $A$ is symmetric.
Assume that $b_1, b_2\in B(\infty)$ satisfy the conditions{\rm:}
\bna
\item one of  $G^\up(b_1)^2$ and $G^\up(b_2)^2$ is a member of the upper global basis
up to a power of $q$,
\item
$G^\up(b_1)$ and $G^\up(b_2)$ q-commute.
\ee
Then their product $G^\up(b_1)G^\up(b_2)$ is a member of the upper global basis 
of $\Um$ up to a power of $q$.
\encor

\medskip
\subsection{Geometric results}\label{subsec:geometry}
The result of this subsection (Theorem~\ref{th:head}) 
is informed us by Peter McNamara.

{\em In this subsection, we assume further that 
the base field $\cor$ is a field of characteristic $0$}.

\Th[{\cite[Lemma 7.5]{Mc14}}]\label{th:head}
Assume that the base field\/ $\cor$ is a field of characteristic $0$.
Assume that $M\in R\gmod$ has
a head $q^cH$ with a self-dual simple module $H$ and $c\in\Z$. Then we have
the equality in the Grothendieck group $K(R\gmod)$
$$[M]=q^c[H]+\sum_{k}q^{c_k}[S_k]$$
with self-dual simple modules $S_k$ and $c_k>c$.
\enth

By duality, we obtain the following corollary.
\Cor Assume that the base field $\cor$ is a field of characteristic $0$.
Assume that $M\in R\gmod$ has a socle $q^cS$
with a self-dual simple module $S$ and $c\in\Z$. Then we have
the equality in $K(R\gmod)$
$$[M]=q^c[S]+\sum_{k}q^{c_k}[S_k]$$
with self-dual simple modules $S_k$ and $c_k<c$.
\encor

Applying this theorem to convolution products, 
we obtain the following corollary.

\Cor Assume that the base field $\cor$ is a field of characteristic $0$.
Let $M$ and $N$ be simple modules.
We assume that one of them is real.
Then we have the equalities in  $K(R\gmod)$:
\bnum
\item
$[M\conv N]=[M\hconv N]+\sum_{k}q^{c_k}[S_k]$\\
with self-dual simple modules $S_k$ and 
$$c_k>-\tL(M,N)=\bl -\La(M,N)-(\wt(M),\wt(N)\br/2.$$
\item
$[M\conv N]=[q^{\La(N,M)}N\hconv M]+\sum_{k}q^{c_k}[S_k]$\\
with self-dual simple modules $S_k$ and $c_k<\bl \La(N,M)-(\wt(N),\wt(M))\br/2$.
\ee
\encor
Note that $q^{\tL(M,N)}M\hconv N$ is self-dual by Lemma~\ref{lem:selfdual}.

Theorem~\ref{th:leclerc} and Theorem~\ref{th:head}
solve affirmatively Conjecture~1
of Leclerc (\cite{L03}) in the symmetric generalized Cartan matrix case.

We obtain the following result
which is a generalization of Lemma~\ref{cor:crde} in
the characteristic-zero case.

\Cor \label{cor:minmaxD}
Assume that the base field $\cor$ is a field of characteristic $0$.
Let $M$ and $N$ be simple modules.
We assume that one of them is real. Write
$$[M\conv N]=\sum_{k=1}^nq^{c_k}[S_k]$$
with self-dual simple modules $S_k$ and $c_k\in\Z$.
Then we have
$$\max\set{c_k}{1\le k\le n}-\min\set{c_k}{1\le k\le n}=\de(M,N).$$
\encor

\subsection{Proof of Theorem~\ref{th:head}}
Recall that the graded algebra $R(\beta)$ ($\beta\in \rtl^+$)
is geometrically realized as follows (\cite{VV09}).
There exist a reductive group $G$ and a $G$-equivariant projective morphism
$f\cl X\to Y$ from a smooth algebraic $G$-variety $X$
to an affine $G$-variety $Y$ defined over the complex number field $\C$
such that
$$R(\beta)\simeq \tEnd_{\Db(\cor_Y)}(\R f_*\cor_X)
\quad\text{as a graded $\cor$-algebra.}$$
Here, $\Db(\cor_Y)$ denotes the equivariant derived category 
of the $G$-variety $Y$ with coefficient $\cor$,
and $\tEnd_{\Db(\cor_Y)}(K)=\tHom_{\Db(\cor_Y)}(K, K)$  with
$$\tHom_{\Db(\cor_Y)}(K, K')\seteq\soplus_{n\in\Z} \Hom_{\Db(\cor_Y)}(K, K'[n]).$$
By the decomposition theorem (\cite{BBDG}), we have a decomposition
$$\R f_*\cor_X\simeq\soplus_{a\in J}E_a\tens \F_a$$
where $\{\F_a\}_{a\in J}$ is a finite family of simple perverse sheaves
on $Y$ and $E_a$ is a finite-dimensional graded $\cor$-vector space such that
\eq
H^k(E_a)\simeq H^{-k}(E_a)\quad\text{for any $k\in\Z$.}
\label{eq:Eselfdual}
\eneq
The last fact \eqref{eq:Eselfdual} follows from the hard Lefschetz theorem
(\cite{BBDG}).

Set $A_{a,b}=\tHom_{\Db(\cor_Y)}(\F_b, \F_a)$.
Then we have the multiplication morphisms
$$A_{a,b}\tens A_{b,c}\to A_{a,c}$$ so that
$$A\seteq\soplus_{a,b\in J}A_{a,b}$$
has a structure of $\Z$-graded algebra such that
$$A_{\le0}\seteq\soplus_{n\le0}A_n=A_0\simeq\cor^J.$$
Hence the family of the isomorphism classes of simple objects
(up to a grading shift) in $A\gmod$
is $\{\cor_a\}_{a\in J}$. 
Here, $\cor_a$ is the module obtained by the algebra homomorphism
$A\to A_{\le0}\simeq \cor^J\to\cor$,
where the last arrow is the $a$-th projection.
Hence we have
$$K(A\gmod)\simeq\soplus_{a\in J}\Z[q^{\pm1}][\cor_a].$$
On the other hand, we have
$$R(\beta)\simeq \soplus_{a,b\in J}E_a\tens A_{a,b}\tens E_b^*.$$
Set
$$L\seteq\soplus_{a,b\in J}E_a\tens A_{a,b}.$$
Then, $L$ is endowed with a natural structure of
$(\soplus_{a,b\in J}E_a\tens A_{a,b}\tens E_b^*, A)$-bimodule.
It is well-known that
the functor $M\mapsto L\otimes_AM$ gives a graded Morita-equivalence
$$\Phi\cl A\gmod\isoto R(\beta)\gmod.$$
Note that
$\Phi(\cor_a)\simeq E_a$ and
$\{E_a\}_{a\in J}$ is the set of isomorphism classes of self-dual
simple graded $R(\beta)$-modules by
\eqref{eq:Eselfdual}.

By the above observation, in order to prove the theorem,
it is enough to show the corresponding statement for
the graded ring $A$, which is obvious.
%
\hfill \qedsymbol

\section{Quantum cluster algebras}
In this section we recall the definition of skew-symmetric quantum cluster algebras following \cite{BZ05}, \cite[\S\,8]{GLS}.

\subsection{Quantum seeds} Fix a finite index set
$\K=\Kex\sqcup\Kfr$ with the decomposition into 
the set of exchangeable indices and the set of frozen indices.
Let $L=(\lambda_{ij})_{i,j\in \K}$ be a skew-symmetric integer-valued $\K\times \K$-matrix.
\Def We define $\mathscr P(L)$ as the  $\Z[q^{\pm 1/2}]$-algebra generated by 
a family of elements $\{X_i\}_{i\in \K}$ 
with the defining relations
\eq
&&X_iX_j=q^{\lambda_{ij}}X_j X_i \quad (i,j\in \K).\label{eq:Xcom}
\eneq
We denote by $\mathscr F(L)$ the skew field of fractions of $\mathscr P(L)$.

\edf
For ${\bf a}=(a_i)_{i\in\K}\in \Z^\K$, we define the element $X^{\bf a}$
of $\mathscr F(L)$
as
\eqn
&&X^{\bf a}\seteq q^{1/2 \sum_{i > j} a_ia_j\lambda_{ij}}
\overset{\To[{\ }]}{\prod}_{i\in\K}
 X_i^{a_i}.
\eneqn
Here we take a total order $<$ on the set $\K$ and
$\overset{\To[{\ }]}{\prod}_{i\in\K} X_i^{a_i}=X_{i_1}^{a_{i_1}}\cdots X_{i_r}^{a_{i_r}}$
where $\K=\{i_1,\ldots,i_r\}$ with $i_1<\cdots <i_r$.
Note that 
$X^{\bf a}$ does not depend on the choice of a total order of $\K$.

We have
\eq
&&X^{\bf a}X^{\bf b}=q^{1/2\sum_{i,j\in\K} a_ib_j\la_{ij}}X^{\bf a+\bf b}.
\label{eq:XaXb}
\eneq
If $\mathbf{a}\in\Z_{\ge0}^\K$, then $X^{\mathbf a}$ belongs to $\mathscr P(L)$.

It is well known that $\{X^{\mathbf a}\}_{\mathbf{a}\in\Z_{\ge0}^\K}$ is a basis  
of $\mathscr P(L)$ as a $\Z[q^{\pm 1/2}]$-module.

Let $A$ be a $\Z[q^{\pm1/2}]$-algebra. We say that a family
$\{x_i\}_{i\in\K}$ of elements of $A$ is {\em $L$-commuting} if
it satisfies $x_ix_j=q^{\la_{ij}}x_jx_i$ for any $i,j\in\K$. 
In such a case we can define $x^\mathbf{a}$ for any $\mathbf{a}\in\Z_{\ge0}^\K$.
We say that an $L$-commuting family $\{x_i\}_{i\in\K}$ is {\em algebraically independent} if the algebra map
$\mathscr P(L)\to A$ given by $X_i\mapsto x_i$ is injective.

Let $\widetilde B = (b_{ij})_{(i,j)\in\K\times\Kex}$ be an  integer-valued 
$\K \times \Kex$-matrix.
 We assume that the {\em principal part $B\seteq(b_{ij})_{i,j\in\Kex}$}  of $\widetilde B$ is skew-symmetric.
We say that the pair $(L, \widetilde B)$ is {\em compatible}, if there exists a positive integer $d$ such that
\eq \label{eq:compatible}
&&\sum_{k\in\K} \lambda_{ik} b_{kj} =\delta_{ij} d \quad (i\in \K,\;j\in\Kex).
\eneq

Let $(L,B)$ be a compatible pair and $A$ a $\Z[q^{\pm1/2}]$-algebra.
We say that $\seed = (\{x_i\}_{i\in\K},L, \widetilde B)$ is
a {\em quantum seed} in $A$ if $\{x_i\}_{i\in\K}$ is 
an algebraically independent $L$-commuting family of elements of
 $A$.

The set $\{x_i\}_{i\in\K}$ is called the {\em cluster} of $\seed$ and
its elements {\em cluster variables}.
The cluster variables $x_i$ ($i\in\Kfr$) are called  {\em frozen variables}.
The elements $x^{\bf a}$ (${\bf a}\in \Z_{\ge0}^\K$) are called {\em quantum cluster monomials}.

\subsection{Mutation}
For $k\in\Kex$, we define a $\K\times \K$-matrix  $E=(e_{ij})_{i,j\in\K}$ and a 
$\Kex\times \Kex$-matrix $F=(f_{ij})_{i,j\in\Kex}$ as follows:
\eqn
e_{ij}=
\begin{cases}
  \delta_{ij} & \text{if} \ j \neq k, \\
  -1 & \text{if} \ i= j = k, \\
  \max(0, -b_{ik}) & \text{if} \ i \neq  j = k,
\end{cases}
\hs{10ex}
f_{ij}=
\begin{cases}
  \delta_{ij} & \text{if} \ i \neq k, \\
  -1 & \text{if} \ i= j = k, \\
  \max(0, b_{kj}) & \text{if} \ i = k \neq j.
\end{cases}
\eneqn
The {\em mutation $\mu_k(L,\widetilde B)\seteq(\mu_k(L),\mu_k(\widetilde B))$ of a compatible pair $(L,\widetilde B)$ in direction $k$} is given by
\eqn
\mu_k(L)\seteq(E^T) \, L \, E, \quad \mu_k(\widetilde B)\seteq E \, \widetilde B \, F.
\eneqn
Then the pair $(\mu_k(L),\mu_k(\widetilde B))$
is also compatible with the same integer $d$ as in the case of $(L,\widetilde B)$ (\cite{BZ05}).

Note that for each $k\in\Kex$, we have
\eq \label{eq:mutatuion B}
\mu_k(\widetilde B)_{ij} =
\begin{cases}
  -b_{ij} & \text{if}  \ i=k \ \text{or} \ j=k, \\
  b_{ij} + (-1)^{\delta(b_{ik} < 0)} \max(b_{ik} b_{kj}, 0) & \text{otherwise,}
\end{cases}
\eneq
and
\eqn
\mu_k(L)_{ij} =
\begin{cases}
  0 & \text{if}  \ i=j \\
  -\la_{kj}+\displaystyle\sum _{t\in\K} \max(0, -b_{tk}) \la_{tj} & \text{if} \ i=k, \ j\neq k, \\
  -\la_{ik}+\displaystyle\sum _{t\in\K} \max(0, -b_{tk}) \la_{it} & \text{if} \ i \neq k, \ j= k, \\
  \la_{ij} & \text{otherwise.}
\end{cases}
\eneqn
Note that we have $$\displaystyle\sum _{t\in\K} \max(0, -b_{tk}) \la_{it}
 =\displaystyle\sum _{t\in\K} \max(0, b_{tk}) \la_{it}$$
 for $i\in\K$ with $i\neq k$, since $(L,\widetilde B)$ is compatible.

We define
\eq
&&a_i'=
\begin{cases}
  -1 & \text{if} \ i=k, \\
 \max(0,b_{ik}) & \text{if} \ i\neq k,
\end{cases} \qquad
a_i''=
\begin{cases}
  -1 & \text{if} \ i=k, \\
 \max(0,-b_{ik}) & \text{if} \ i\neq k.
\end{cases}
\label{eq:aa}
\eneq
and set ${\bf a}'\seteq(a_i')_{i\in\K}$ and ${\bf a}''\seteq(a_i'')_{i\in\K}$.

Let $A$ be a $\Z[q^{\pm1/2}]$-algebra contained in a skew field $K$.
Let $\seed=(\{x_i\}_{i\in\K},L, \widetilde B)$ be a quantum seed in $A$.
Define the elements $\mu_k(x)_i$ of $K$ by
\eq \label{eq:quantum mutation}
\mu_k(x)_i\seteq
\begin{cases}
  x^{{\bf a}'}  +   x^{{\bf a}''}, & \text{if} \ i=k, \\
 x_i & \text{if} \ i\neq k.
\end{cases}
\eneq
Then $\{\mu_k(x)_i\}$ is an algebraically independent
$\mu_k(L)$-commuting family in $K$.
We call 
\eqn
\mu_k(\seed)\seteq\bl\{\mu_k(x)_i\}_{i\in\K},\mu_k(L),\mu_k(\widetilde B)\br
\eneqn
the {\em mutation 
of $\seed$ in direction $k$}.
It becomes a new quantum seed in $K$.

 \begin{definition}
Let $\seed=(\{x_i\}_{i\in\K},L, \widetilde B)$ be a quantum seed in $A$.
   The {\em quantum cluster algebra $\mathscr A_{q^{1/2}}(\seed)$} associated to the quantum seed $\seed$ is  the $\Z[q^{\pm 1/2}]$-subalgebra of the skew field $K$ generated by all the quantum cluster variables in the quantum seeds obtained from $\seed$ by any sequence of mutations.
 \end{definition}
We call $\seed$ the {\em initial quantum seed} 
of the quantum cluster algebra $\mathscr A_{q^{1/2}}(\seed)$.

\section{Monoidal categorification}

Throughout  this section, fix $\K=\Kex\sqcup\Kfr$, and a base field
$\coro$.

Let $\shc$ be a $\coro$-linear abelian monoidal category. 
For the definition of monoidal category, see \cite[Appendix A.1]{K^3}. 
Note that  in \cite{K^3}, it was called  \emph{tensor category}.  
A  $\coro$-linear abelian monoidal category is a $\coro$-linear 
monoidal category such that
it is abelian and
the tensor functor $\;\tens\;$ is $\coro$-bilinear and exact.

We assume further the following conditions on $\shc$ 
\eq
&&\hs{-4ex}\left\{\parbox{60ex}{
\bnum
\item
Any object of $\shc$ has a finite length,
\vs{.5ex}
\item $\coro\isoto\Hom_{\shc}(S,S)$ for any simple object $S$ of $\shc$.
\ee}\right.\label{cond:monoidal category}
\eneq

\subsection{Ungraded cases}
\begin{definition}
  \label{def:monoidal seed}
Let  $\seed = (\{ M_i\}_{i\in \K },\widetilde B)$ be a pair
of 
a family $\{ M_i\}_{i\in\K}$ of simple objects in $\shc$ and
an  integer-valued $\K\times\Kex$-matrix 
$\widetilde B = (b_{ij})_{(i,j)\in\K\times\Kex}$
whose principal part is skew-symmetric.
We call $\seed$ a \emph{monoidal seed in $\shc$} if
\bnum
  \item $M_i \otimes M_j \simeq M_j \otimes M_i$ for any $i,j\in\K$, and
  \item $\sotimes_{i\in\K} M_i^{\otimes a_i}$ is simple for any $(a_i)_{i\in\K}\in \Z_{\ge 0}^{\K}$.
  \end{enumerate}
\end{definition}

\begin{definition}
  \label{def:mutation monoidal seed}
  For $k\in\Kex$,
  we say that a  monoidal seed  $\seed = (\{ M_i\}_{i\in \K },\widetilde B)$
 \emph{admits a mutation in direction $k$} if
there exists a simple object  $M_k' \in \shc$
such that
\bna
  \item
there exist exact sequences in $\shc$
\eq
&&0 \to  \sotimes_{b_{ik} >0} M_i^{\tensor b_{ik}} \to M_k \tensor M_k' \to
 \sotimes_{b_{ik} <0} M_i^{\tensor (-b_{ik})} \to 0, \label{eq:ses_mutation1}\\
 &&0 \to  \sotimes_{b_{ik} <0} M_i^{\tensor(-b_{ik})} \to M_k' \tensor M_k \to
  \sotimes_{b_{ik} >0} M_i^{\tensor b_{ik}} \to 0.\label{eq:ses_mutation2}
\eneq
\item the pair $\mu_k(\seed)\seteq
(\{M_i\}_{i\neq k}\cup\{M_k'\},\mu_k(\widetilde B))$ is
a monoidal seed in $\shc$.
\end{enumerate}
\end{definition}

Recall that a cluster algebra $A$ with an initial seed 
$(\{x_i\}_{i\in\K},\widetilde B)$ is the $\Z$-subalgebra of 
$\Q(x_i\vert i\in\K)$ generated by all the cluster variables in the seeds obtained from  $(\{x_i\}_{i\in\K},\widetilde B)$ by any sequence of mutations.
Here, the mutation $x_k'$ of a cluster variable $x_k$ $(k\in\Kex)$ is given  by
\eqn
x_k' = \dfrac{\prod_{b_{ik} \ge 0} x_i^{b_{ik}}+\prod_{b_{ik} \le 0} x_i^{-b_{ik}}}{x_k},
\eneqn
and the mutation of $\widetilde B$ is given in \eqref{eq:mutatuion B}.

\begin{definition}
A $\coro$-linear abelian monoidal category $\shc$ 
with \eqref{cond:monoidal category} is called a \em{monoidal categorification of a cluster algebra $A$}
if
\bna
\item the Grothendieck ring $K(\shc)$ is isomorphic to $A$, 
\item there exists a monoidal seed $\seed = ( \{M_i\}_{i\in\K},\widetilde B)$ in $\shc$ such that
$[\seed]\seteq( \{[M_i]\}_{i\in\K},\widetilde B)$ is the initial seed of $A$ and 
$\seed$ admits  successive mutations in all directions.
\end{enumerate}
\end{definition}
Note that if $\shc$ is a monoidal categorification of $A$, 
then every seed in $A$ is of the form
$( \{[M_i]\}_{i\in\K},\widetilde B)$ for some monoidal seed 
$( \{M_i\}_{i\in\K},\widetilde B)$.
 In particular, all the cluster monomials in $A$ are the classes of real
simple objects in $\shc$.

\subsection{Graded cases}
Let $\rootl$ be a free abelian group equipped with a symmetric bilinear form
\eqn
( \ , \ ) : \rootl \times \rootl \to \Z \quad \text{such that} \
 (\beta,\beta) \in 2\Z \ \text{for all} \ \beta \in \rootl.
\eneqn
We consider  a $\coro$-linear abelian monoidal category $\shc$ satisfying 
\eqref{cond:monoidal category} and the following conditions:
\eq
&&\hs{-2ex}\left\{\parbox{67ex}{
\bnum
\item We have a direct sum decomposition 
  $\shc = \soplus_{\beta \in \rootl} \shc_\beta $ such that
 the tensor product $\otimes$ 
   sends  $\shc_\beta \times \shc_\gamma$ to $\shc_{\beta + \gamma}$ 
   for every  $\beta, \gamma \in \rootl$.

 \item There exists an object $Q  \in \shc_0$ satisfying
\bna
\item there is an isomorphism
$$R_{Q}(X) : Q \tensor X \isoto X \tensor Q$$
functorial in $X \in \shc$ such that
$$\xymatrix@C=7ex{
Q\tensor X\tensor Y\ar[r]_-{R_Q(X)}\ar@/^3ex/[rr]^{R_Q(X\tensor Y)}
&X\tensor Q\tensor Y\ar[r]_-{R_Q(Y)}
&X\tensor Y\tensor Q}$$
 commutes for any $X,Y\in \shc$,
\item the functor $X \mapsto Q \tensor X$ is an  equivalence of categories.
\end{enumerate}
\item for any $M$, $N\in\shc$, we have
$\Hom_\shc(M,Q^{\tens n}\tens N)=0$ except finitely many integers $n$.
\end{enumerate}
}\right.\label{cond:quantum monoidal category}
\eneq

We denote by $q$ the auto-equivalence $Q \otimes \scbul$ of $\shc$,
and call it the {\em grading shift functor}.

In such a case the Grothendieck group $K(\shc)$ is a $\rtl$-graded
$\Z[q^{\pm1}]$-algebra: $K(\shc)=\soplus_{\beta\in\rtl}K(\shc)_\beta$ where
$K(\shc)_\beta=K(\shc_\beta)$. Moreover we have
$$K(\shc)=\soplus_{S}\Z[q^{\pm1}][S]$$
where $S$ ranges over the set of equivalence classes of
simple modules $S$. Here, two simple modules $S$ and $S'$ are equivalent
if $q^nS\simeq S'$ for some $n\in\Z$.

\medskip
For $M\in\shc_{\beta}$, we sometimes write
$\beta=\wt(M)$ and call it the {\em weight} of $M$.
Similarly, for $x\in \Q(q^{1/2})\tens_{\Z[q^{\pm1}]}K(\shc_\beta)$,
we write $\beta=\wt(x)$ and call it the {\em weight} of $x$.

\begin{definition} \label{def:quantum monoidal seed}
We call a quadruple $\seed = (\{M_i\}_{i\in\K}, L,\widetilde B, D)$
a \emph{quantum monoidal seed  in $\shc$}
 if it satisfies the following conditions{\rm:}
\bnum
\item $\widetilde B = (b_{ij})_{i\in\K,\,j\in\Kex}$ is an  integer-valued $\K\times\Kex$-matrix 
whose principal part is skew-symmetric,
\item $L=(\lambda_{ij\in\K})$ is an integer-valued skew-symmetric $\K\times\K$-matrix,
\item $D=\{d_i\}_{i\in\K}$ is a family of elements in $\rootl$,
\item $\{M_i\}_{i\in\K}$ is a family of simple objects such that $M_i \in \shc_{d_i}$ for any $i\in\K$,
\item $M_i \otimes M_j \simeq q^{\lambda_{ij}} M_j \otimes M_i$ for all $ i, j \in\K$,
\item $M_{i_1} \otimes \cdots \otimes M_{i_t}$ is simple for any
sequence  $(i_1,\ldots,i_t)$ in $\K$,
\item The pair  $(L,\widetilde B)$  is  compatible in the sense of \eqref{eq:compatible} with $d=2$,
\item $\lambda_{ij} - (d_i,d_j) \in 2\Z$ for all $i,j \in\K$,
\item $\displaystyle\sum_{i\in\K}b_{ik}d_i =0$ for all $k\in\Kex$.
\end{enumerate}
\end{definition}

Let  $\seed =(\{M_i\}_{i\in\K}, L,\widetilde B, D)$
be a  quantum monoidal seed.
For any $X\in\shc_\beta$ and $Y\in\shc_\gamma$ such that $X\conv Y\simeq
q^{c}Y\circ X$ and $c+(\beta,\gamma)\in2\Z$, we set
\eq&&\tL(X,Y)=\frac{1}{2}\bl-c+(\beta,\gamma)\br\in\Z\eneq
and
\eq X\nconv Y\seteq q^{\tL(X,Y)}X\tens Y\simeq
q^{\tL(Y,X)}Y\tens X.\label{eq:balpro}\eneq
Then $X\nconv Y\simeq Y\nconv X$.
For any sequence $(i_1,\ldots, i_\ell)$ in $\K$,
we define
$$\sodot_{s=1}^{\ell} M_{i_s}\seteq
(\cdots((M_{i_1}\nconv M_{i_2})\nconv M_{i_3})\cdots)\nconv M_{i_\ell}.$$
Then we have
\eqn
\sodot_{s=1}^{\ell} M_{i_s} =
q^{\frac{1}{2}\sum_{1\le u<v \le \ell}(-\lambda_{i_u i_v}  +(d_{i_u},d_{i_v}) ) } M_{i_1} \tensor \cdots \tensor M_{i_\ell}.
\eneqn
For any $w\in\sym_\ell$, we have
\eqn
&&\sodot_{s=1}^{\ell} M_{i_w(s)}\simeq
\sodot_{s=1}^{\ell} M_{i_s}\eneqn
Hence for any subset $A$ of $\K$ and any set of non-negative integers
$\{m_a\}_{a\in A}$, we can define
$\sodot_{a\in A}M_a^{\snconv m_a}$.

For $(a_i)_{i\in\K}\in\Z_{\ge0}^\K$ and $(b_i)_{i\in\K}\in\Z_{\ge0}^\K$, 
we have
$$\bl\sodot_{i\in \K}M_i^{\snconv a_i}\br\nconv
\bl\sodot_{i\in\K}M_i^{\snconv b_i}\br\simeq\sodot_{i\in \K}M_i^{\snconv (a_i+b_i)}.$$

\medskip
Let $\seed=(\{M_i\}_{i\in\K}, L,\widetilde B, D)$ be  a quantum monoidal seed.
When the $L$-commuting family $\{[M_i]\}_{i\in\K}$ 
of elements of $\Z[q^{\pm1/2}]\tens_{\Z[q^{\pm1}]}K(\shc)$ 
is algebraically independent,
we shall define a quantum seed
$[\seed]$ in $\Z[q^{\pm1/2}]\tens_{\Z[q^{\pm1}]}K(\shc)$
by
\eq&&[\seed]=(\{q^{-(d_i,d_i)/4}[M_i]\}_{i\in\K}, L,\widetilde B).
\eneq

Set 
\eq&&X_i=q^{-(d_i,d_i)/4}[M_i].\label{eq:XM}
\eneq
Then for any $\mathbf a=(a_i)_{i\in\K}\in\Z_{\ge0}^\K$, we have
$$[\sodot_{i\in \K}M_i^{\snconv a_i}]=q^{(\mu,\mu)/4}X^{\mathbf a}$$
where $\mu=\wt(\sodot_{i\in \K}M_i^{\snconv a_i})=\wt(X^{\mathbf a})=\sum_{i\in K}a_id_i$.

\medskip
For a given $k\in\Kex$, we define the \emph{mutation $\mu_k(D) \in \rootl^\K$ of $D$ in direction $k$ with respect to $\widetilde B$} by
\eqn
\mu_k(D)_i =d_i \ (i \neq k), \quad \mu_k(D)_k=-d_k+\sum_{b_{ik} >0}   b_{ik} d_i.
\eneqn
Note that
\eqn
\mu_k(\mu_k(D))=D.
\eneqn
Note also  that 
$(\mu_k(L),\mu_k(\tB),\mu_k(D))$ satisfies conditions (viii) and (ix)
in Definition~\ref{def:quantum monoidal seed} for any $k\in\Kex$.

\medskip

Note that
\Lemma\label{lem:decat} Set $X'_k=\mu_k(X)_k$, the mutation of $X_k$ 
as in \eqref{eq:quantum mutation}. 
Set $\zeta=\wt(X'_k)=-d_k+\sum_{b_{ik}>0}b_{ik}d_i$.
Then we have
\eq
&&\hs{-20ex}\ba{l}q^{m_k}[M_k] q^{(\zeta,\zeta)/4}X'_k=
q [\sodot_{b_{ik} >0} M_i^{\snconv b_{ik}}]+
[\sodot_{b_{ik} <0} M_i^{\snconv (-b_{ik})}],\\[2.5ex]
q^{m'_k} q^{(\zeta,\zeta)/4}X'_k[M_k]= [\sodot_{b_{ik} >0} M_i^{\snconv b_{ik}}]+
q[\sodot_{b_{ik} <0} M_i^{\snconv (-b_{ik})}],\ea
\eneq
where
\eq
&&\hs{-10ex}\left\{\ba{l} 
m_k=\dfrac{1}{2}(d_k,\zeta) +\dfrac{1}{2} \displaystyle 
\sum_{b_{ik} < 0}\la_{ki}  b_{ik},
\\[1ex]
m'_k=\dfrac{1}{2}(d_k,\zeta) +\dfrac{1}{2} \displaystyle 
\sum_{b_{ik} > 0} \la_{ki}b_{ik}. 
\ea\right.\label{eq:mm'}
\eneq
\enlemma
\Proof
By \eqref{eq:XaXb}, we have 
\eqn
X_k X^{\bf a} = q^{\frac{1}{2} \sum_{i \in \K} a_i \la_{ki}}  X^{{{\bf e}_k} +\bf a} 
\quad \text{for ${\bf a}=(a_i)_{i\in\K} \in \Z^\K$ and 
$({\bf e}_k)_i = \delta_{ik}$  ($i \in \K$).}
\eneqn
Let $\mathbf{a}'$ and ${\bf a}''$ be as in \eqref{eq:aa}. Because
\eqn
\sum_{i \in \K} a_i' \la_{ki} - \sum_{i \in \K} a_i'' \la_{ki} 
= \sum_{b_{ik} >0} b_{ik} \la_{ki} - \sum_{b_{ik} <0} (-b_{ik}) \la_{ki} 
= \sum_{i \in \K} b_{ik} \la_{ki} = 2,
\eneqn
 we have
\eqn 
X_k X'_k = X_k(X^{{\bf a}'}  +   X^{{\bf a}''})
=q^{\frac{1}{2} \sum_i a_i'' \la_{ki} } (qX^{{\bf e_k}+{\bf a}'} + X^{{\bf e_k}+{\bf a}''}).
\eneqn
Note that 
$
\wt(X^{{\bf e_k}+{\bf a}'}) =\wt(X^{{\bf e_k}+{\bf a}''})= d_k +\zeta.
$
It follows that 
\eqn
m_k&=&-\dfrac{1}{4}((d_k,d_k)+(\zeta,\zeta))-\dfrac{1}{2} \sum_{i \in J} a_i'' \la_{ki} + \dfrac{1}{4}(\zeta+d_k,\zeta+d_k) \\
&=&\dfrac{1}{2} (d_k, \zeta) +\dfrac{1}{2} \sum_{b_{ik} <0} b_{ik} \la_{ki}.
\eneqn
One can calculate $m_k'$ in a similar way.
\QED

\begin{definition} \label{def:monoidal mutation}
We say that a quantum monoidal seed  
$\seed =(\{M_i\}_{i\in\K}, L,\widetilde B, D)$
 \emph{admits a mutation in direction $k\in\Kex$} if
there exists  a simple object  $M_k' \in \shc_{\mu_k(D)_k}$ 
such that
\bna
  \item
there exist exact sequences in $\shc$
\eq
&&0 \to q \sodot_{b_{ik} >0} M_i^{\snconv b_{ik}} \to q^{m_k} M_k \tensor M_k' \to
 \sodot_{b_{ik} <0} M_i^{\snconv (-b_{ik})} \to 0, 
\label{eq:ses graded mutation1} \\
 &&0 \to q \sodot_{b_{ik} <0} M_i^{\snconv(-b_{ik})} \to q^{m_k'} M_k' \tensor M_k \to
  \sodot_{b_{ik} >0} M_i^{\snconv b_{ik}} \to 0, \label{eq:ses graded mutation2}
\eneq
where $m_k$ and $m'_k$ are as in \eqref{eq:mm'}.
\item 
$\mu_k(\seed)\seteq\bl\{M_i\}_{i\neq k}\sqcup\{M_k'\},\mu_k(L),
\mu_k(\widetilde B), \mu_k(D)\br$ is
a quantum monoidal seed in $\shc$.
\end{enumerate}
We call $\mu_k(\seed)$ the {\em mutation} of $\seed$ in direction $k$.
\end{definition}

By Lemma~\ref{lem:decat}, the following lemma is obvious.
\Lemma
Let $\seed =(\{M_i\}_{i\in\K}, L,\widetilde B, D)$ be
 a quantum monoidal seed  which admits a mutation in direction $k\in\Kex$.
Then we have
$$[\mu_k(\seed)]=\mu_k([\seed]).$$
\enlemma

\begin{definition}
Assume that a $\coro$-linear abelian monoidal category 
$\shc$ satisfies  \eqref{cond:monoidal category}
and \eqref{cond:quantum monoidal category}.
The category $\shc$ is called a 
\em{monoidal categorification of a quantum cluster algebra $A$ 
over $\Z[q^{\pm1/2}]$}
if
\bnum
\item the Grothendieck ring $\Z[q^{\pm1/2}]\tens_{\Z[q^{\pm1}]} K(\shc)$ is isomorphic to $A$, 
\item there exists a quantum monoidal seed 
$\seed =(\{M_i\}_{i\in\K}, L,\widetilde B, D)$ in $\shc$ such that
$[\seed]\seteq(\{q^{-(d_i,d_i)/4}[M_i]\}_{i\in\K}, L, \widetilde B)$
 is a quantum seed of $A$,
\item $\seed$ admits successive mutations in all the directions.
\end{enumerate}
\end{definition}
Note that if $\shc$ is a monoidal categorification of a quantum cluster algebra
$A$, then any quantum seed in $A$ obtained by mutations from the initial quantum seed
is of the form
$(\{q^{-(d_i,d_i)/4}[M_i]\}_{i\in\K}, L, \widetilde B)$  for some quantum monoidal seed 
$(\{M_i\}_{i\in\K}, L,\widetilde B, D)$.
 In particular, all the quantum cluster monomials in $A$ are the
classes of real simple objects in $\shc$ up to a power of $q^{1/2}$.

\section{Monoidal categorification via modules over KLR algebras}
Let $R$ be a symmetric KLR algebra 
over a base field $\coro$.

From now on, we focus on the case when $\shc$ is a full subcategory of $R \gmod$
 stable under taking convolution products, subquotients, extensions
and grading shift.
In particular, we have
\eqn
\shc = \soplus_{\beta \in \rootl^-} \shc_\beta, \quad
\text{where} \ \shc_\beta \seteq\shc\cap R(-\beta) \gmod,
\eneqn
and we have the grading shift functor $q$ on $\shc$.
Hence we have
$$K(\shc_\beta)\subset\Um_\beta,$$
and $K(\shc)$ has a $\Z[q^{\pm1}]$-basis consisting of 
the isomorphism classes of self-dual simple modules.

\begin{definition} \label{def:admissible}
A pair $(\{M_i\}_{i\in\K}, \widetilde B)$ is called \emph{admissible} if
\bna
\item  $\{M_i\}_{i\in\K}$ is a family of real simple  self-dual  objects of $\shc$ which commute with each other,
\item $\widetilde B$ is an integer-valued $\K \times\Kex$-matrix with 
skew-symmetric principal part,
\item
 for each $k\in\Kex$, there exists a  self-dual  simple object $M'_k$ of $\shc$ 
 such that
 there is an exact sequence in $\shc$
 \eq 
&&0 \to q \sodot_{b_{ik} >0} M_i^{\snconv b_{ik}} \to q^{\tLa(M_k,M_k')} M_k \conv M_k' \to
 \sodot_{b_{ik} <0} M_i^{\snconv (-b_{ik})} \to 0, 
\label{eq:ses graded mutation KLR}
 \eneq
and
$M_k'$ commutes with $M_i$  for any  $i \neq k$. 
  \end{enumerate}
\end{definition}

Note that $M'_k$ is uniquely determined by $k$ and 
$(\{M_i\}_{i\in\K}, \widetilde B)$. Indeed, it follows from
$q^{\tLa(M_k,M_k')}M_k\hconv M'_k\simeq\sodot_{b_{ik} <0} M_i^{\snconv(-b_{ik})}$ and
\cite[Corollary 3.7]{KKKO14}.
Note also that if there is an epimorphism $q^mM_k\conv M'_k\epito
\sodot_{b_{ik} <0} M_i^{\snconv(-b_{ik})}$ for some $m\in\Z$, then $m$ should coincide with $\tL(M_k,M'_k)$ by Lemma~\ref{lem:selfdual} and Lemma~\ref{lem:multipro}.

\smallskip
For an admissible pair  $(\{M_i\}_{i\in\K}, \widetilde B)$, let
$\La=(\La_{ij})_{i,j\in\K}$
be the skew-symmetric matrix
given by $\La_{ij}=\Lambda(M_i,M_j)$.
and let $D=\{d_i\}_{i\in\K}$ be the family of elements of $\rootl^-$ given by 
$d_i=\wt(M_i)$.

\smallskip
Now we can simplify  the conditions in Definition \ref{def:quantum monoidal seed} and
Definition \ref{def:monoidal mutation} as follows.

\begin{prop} \label{prop:condition simplified}
Let $(\{M_i\}_{i\in\K},\widetilde B)$ be an admissible pair in $\shc$,
and let $M'_k$ $(k\in\Kex)$ be as in {\rm Definition~\ref{def:admissible}}.
    Then we have the following properties.
  \bnum
      \item The quadruple $\seed\seteq(\{M_i\}_{i\in\K}, -\La,\widetilde B,D)$
    is a quantum monoidal seed in $\shc$.\label{item:1}
    \item The self-dual simple object $M_k'$ is real for every $k\in\Kex$.
\label{item:2}
    \item The quantum monoidal seed $\seed$ admits a mutation in each direction $k\in\Kex$.\label{item:3}
\item $M_k$ and $M'_k$ is simply-linked for any $k\in\Kex$
\ro i.e., $\de(M_k,M'_k)=1$\rf.\label{item:sl}
\item For any $j\in\K$ and $k\in\Kex$, we have
\eq&&
\ba{l}\La(M_j,M'_k)=-\La(M_j,M_k)-\sum_{b_{ik}<0}\La(M_j,M_i)b_{ik},\\[1ex]
\La(M'_k,M_j)=-\La(M_k,M_j)+\sum_{b_{ik}>0}\La(M_i,M_j)b_{ik}.
\ea\label{eq:Larel}
\eneq\label{item:Larel}
  \end{enumerate}
\end{prop}
\Proof
\eqref{item:sl} follows from the exact sequence
\eqref{eq:ses graded mutation KLR} and Corollary \ref{cor:crde}.

\medskip\noi
\eqref{item:2} follows from the exact sequence
\eqref{eq:ses graded mutation KLR}
by applying Corollary \ref{cor:real} to the case
$$M=M_k, \ N=M'_k, \ X= q  \sodot_{b_{ik} > 0} M_i^{\snconv b_{ik}},
\ \text{and} \  Y= \sodot_{b_{ik} < 0} M_i^{\snconv(-b_{ik})}. $$

\medskip
\noi
\eqref{item:Larel} follows from
\begin{align*}
\La(M_j,M_k)+\La(M_j,M'_k)&=\La(M_j,M_k\hconv M'_k)
=\La\bl M_j,\sodot_{b_{ik} <0} M_i^{\snconv(-b_{ik})}\br\\
&=\sum_{b_{ik}<0}\La(M_j,M_i)(-b_{ik})
\end{align*}
and
\begin{align*}
\La(M_k,M_j)+\La(M'_k,M_j)&=\La(M'_k\hconv M_k,M_j)
=\La\bl\sodot_{b_{ik} >0} M_i^{\snconv b_{ik}}, M_j\br\\
&=\sum_{b_{ik}>0}\La(M_i,M_j)b_{ik}.
\end{align*}

\medskip\noi
Let us show \eqref{item:1}.
Conditions (i)--(v) in Definition \ref{def:quantum monoidal seed} are satisfied by the construction.
Condition (vi) follows from Proposition \ref{prop:real simple}
and the fact that $M_i$ is real simple 
for every $i\in\K$.
Condition (viii) is nothing but Lemma \ref{lem:tLa even}.
Condition (ix) follows easily from the fact that
the weights of the first and the last terms in the exact sequence \eqref{eq:ses graded mutation KLR} coincide.

\medskip\noi
Let us show condition (vii)  in Definition \ref{def:quantum monoidal seed}.
By \eqref{eq:Larel} and \eqref{item:sl} of this proposition, we have
\begin{align*}
2\delta_{jk}=2\de(M_j,M'_k)
&=-2\de(M_j,M_k)-\sum_{b_{ik}<0}\La(M_j,M_i)b_{ik}
+\sum_{b_{ik}>0}\La(M_i,M_j)b_{ik}\\
&=-\sum_{b_{ik}<0}\La(M_j,M_i)b_{ik}-\sum_{b_{ik}>0}\La(M_j,M_i)b_{ik}
=-\sum_{i\in\K}\La(M_j,M_i)b_{ik}.
\end{align*}
for $k\in\Kex$ and $j\in\K$.

Thus we have shown that $\seed$ is a quantum monoidal seed 
in $\shc$.

\medskip\noi
Let us show \eqref{item:3}. Let $k\in\Kex$. The exact sequence 
\eqref{eq:ses graded mutation1}  
follows from  \eqref{eq:ses graded mutation KLR}
and the equality
\eq
\tL(M_k,M'_k)=\dfrac{1}{2}\bl(\wt(M_k,M'_k)
-\sum_{b_{ik}<0}\La(M_k,M_i)b_{ik}\br=m_k
\eneq
which is an immediate consequence of \eqref{eq:Larel}.

 Similarly, taking the dual of
the exact sequence \eqref{eq:ses graded mutation KLR}, 
we obtain an exact sequence
 \eqn
&&0 \to  \sodot_{b_{ik} <0} M_i^{\snconv (-b_{ik})}  \to q^{-\tLa(M_k,M_k')+(\wt M_k, \wt M_k')} M_k' \conv M_k  \to q^{-1}\sodot_{b_{ik} >0} M_i^{\snconv b_{ik}}
 \to 0,
 \eneqn
which gives the exact sequence \eqref{eq:ses graded mutation2}.
\smallskip

It remains to prove that $\mu_k(\seed)\seteq(\{M_i\}_{i\neq k}\cup\{M_k'\},\mu_k(-\La),
\mu_k(\widetilde B), \mu_k(D))$ is
a quantum monoidal seed in $\shc$ for any $k\in\Kex$.

We see easily that $\mu_k(\seed)$
satisfies the conditions (i)--(iv) and (vii)--(ix) in Definition~\ref{def:quantum monoidal seed}.

For condition (v), it is enough to show that
for $i,j\in\K$ we have
\eqn \mu_k(-\La)_{ij} = -\Lambda(\mu_k(M)_i,\mu_k(M)_j),
\eneqn
where $\mu_k(M)_i = M_i$ for $i\neq k$ and $\mu_k(M)_k=M_k'$.
In the case $i \neq k$ and $j \neq k$, we have
\eqn \mu_{k}(-\La)_{ij}=-\Lambda(M_i,M_j)=-\Lambda(\mu_k(M_i),\mu_k(M_j)). 
\eneqn
The other cases follow from \eqref{eq:Larel}.

Condition (vi) in Definition~\ref{def:quantum monoidal seed} 
for $\mu_k(\seed)$ follows from Proposition \ref{prop:real simple} and the fact that $\{\mu_k(M)_i\}_{i\in\K}$ is a commuting family of real simple modules.
\QED

Now we are ready to give  our main theorem.
\Th\label{th:main}
Let  $(\{M_i\}_{i\in K},\widetilde B)$ be an admissible pair in $\shc$
and set $$\seed=(\{M_i\}_{i\in\K}, -\La,\widetilde B,D)$$
as in {\rm Proposition~\ref{prop:condition simplified}}.
We set
$[\seed]\seteq(\{q^{-\frac{1}{4}(\wt(M_i), \wt(M_i))}[M_i]\}_{i\in\K}, -\La,\widetilde B,D)$.
We assume further
\eq
&&\hs{1ex}\parbox[t]{79ex}{
 The $\Q(q^{1/2})$-algebra $\Q(q^{1/2})\tens\limits_{\Z[q^{\pm1}]}K(\shc)$ 
is isomorphic to $\Q(q^{1/2})\tens\limits_{\Z[q^{\pm1}]}\mathscr A_{q^{1/2}}([\seed])$.
}\label{eq:Cluster}
\eneq
Then, for each $x\in\Kex$, the pair 
$\bl\{\mu_x(M)_i\}_{i\in\K},\mu_x(\widetilde B)\br$
 is admissible in $\shc$.
\enth
\Proof In Proposition ~\ref{prop:condition simplified}, we have already proved 
the conditions (a) and (b) in Definition~\ref{def:admissible} for 
$(\{\mu_x(M)_i\}_{i\in\K},\mu_x(\widetilde B))$.
Let us show (c).
Set $N_i\seteq\mu_x(M)_i$ and
$b_{ij}' \seteq\mu_x(\widetilde B)_{ij}$ for $i\in \K$ and $j\in \Kex$.
It is enough to show that, for any $y\in\Kex$,
there exists a self-dual simple module $M_y'' \in \shc$ such that
there is a short exact sequence
\eq
\xymatrix{
0 \ar[r] &  q \sodot_{b'_{iy} > 0} N_i^{\snconv b'_{iy}} \ar[r]  
&q^{\tLa(N_y,M_y'')}  N_y \conv  M_y'' \ar[r] 
& \sodot_{b'_{iy} <0} N_i^{\snconv (-b'_{iy})}   \ar[r] & 0,
}\label{eq:seqdes}
\eneq
and
we have
$$\de(N_i,M_y'')=0 \quad \text{for $i \neq y$.}$$

If $x=y$, then $b'_{iy}=-b_{ix}$ and hence $M''_y=M_x$ satisfies the desired condition.

Assume that $x\not=y$ and $b_{xy}=0$. Then $b'_{iy} = b_{iy}$ for any $i$
and $N_i=M_i$ for any $i\not=x$.
Hence  $M_y''=\mu_y(M)_y$ satisfies the desired condition.

\medskip
We will show the assertion in the case $b_{xy} > 0$.
We omit the proof of the case $b_{xy} <0$ because it
can be shown in a similar way.

Set
\eqn
&&M_x'\seteq\mu_x(M)_x, \hs{3ex} M_y'\seteq\mu_y(M)_y, \\
&&C\seteq\sodot_{b_{ix} >0} M_i^{\snconv b_{ix}}, \quad  S\seteq\sodot_{b_{ix} <0, \ i\neq y} M_i^{\snconv -b_{ix}}, \\
&&P\seteq\sodot_{b_{iy} >0, i\neq x} M_i^{\snconv b_{iy}}, \quad  Q=\sodot_{b_{iy}' <0, \ i\neq x} M_i^{\snconv-b_{yi}'}, \\
&&A\seteq\sodot_{b_{iy}' \le 0, \ b_{ix} >0} M_i^{\snconv b_{ix}b_{xy}} \nconv \sodot_{\substack{b_{iy} <0,\; b_{iy}' >0,\;b_{ix} >0}} M_i^{\snconv -b_{iy}}
 \\
&&\hs{3ex}\simeq
\sodot_{\substack{b_{iy} <0,\; b_{ix} >0}} 
M_i^{\snconv \min(b_{ix}b_{xy},-b_{iy})},\\
&&B\seteq\sodot_{b_{iy} > 0, \ b_{ix} >0} M_i^{\snconv b_{ix}b_{xy}} \nconv \sodot_{\substack{b_{iy}' >0, \; b_{iy} <0,\;b_{ix} >0}} M_i^{\snconv b_{iy}'}.
 \eneqn

Then we have
\begin{align*}
Q \nconv A &\simeq \sodot_{b_{iy} < 0} M_i^{\snconv -b_{iy}} \quad \text{and}\\
A \nconv B &\simeq C^{\snconv b_{xy}}.
\end{align*}

Set
 \eqn
&&L\seteq(M_x')^{\snconv b_{xy}} \quad \text{and} \quad
V\seteq M_x^{\snconv b_{xy}}.
\eneqn

Set
\eqn
X \seteq \sodot_{b_{iy} >0} M_i^{\snconv b_{iy}}\simeq
M_x^{\snconv b_{xy}} \nconv P =V \nconv P,
 \quad Y\seteq \sodot_{b_{iy} <0} M_i^{\snconv-b_{iy}}\simeq Q \nconv A.
\eneqn

Then \eqref{eq:seqdes} reads as
\eq
\xymatrix{
0 \ar[r] &  q (B \nconv P)\ar[r]  &q^{\tL(M_y,M''_y)}M_y \conv M''_y \ar[r] &
L  \nconv Q\ar[r] & 0.
}\label{eq:seqdes'}
\eneq

Note that we have
\eq
&&0 \to q C  \to q^{\tLa(M_x,M_x')}M_x \conv M_x' \to  M_y^{\snconv b_{xy}} \nconv S \to 0 \quad \text{and} \\
&&0 \to q X \to  q^{\tLa(M_y,M_y')}M_y \conv M_y' \to  Y \to 0. \label{eq:sesMy}
\eneq

Taking the convolution products of $L=(M_x')^{\snconv b_{xy}}$  and \eqref{eq:sesMy}, we obtain
\eqn
&&
\xymatrix@R=.5ex{0 \ar[r]& qL \conv X  \ar[r]&
 q^{\tLa(M_y,M_y')} L \conv (M_y \conv M_y')\ar[r]& L \conv Y \ar[r]& 0, \\
0\ar[r]& qX \conv L \ar[r]& 
q^{\tLa(M_y,M_y')} (M_y \conv M_y') \conv L \ar[r]&  Y\conv L \ar[r]& 0.}
\eneqn

Since $L$ commutes with  $M_y$,
we have
\eqn
&&\Lambda(L, Y)=\Lambda(L, M_y \hconv M_y') \\
 &&=\Lambda(L, M_y) +\Lambda(L, M_y') = \Lambda(L,M_y \conv M_y' ).
\eneqn

On the other hand, we have
\eqn
&&\Lambda(M_x',X) - \Lambda(M_x',Y)  \\
&&=\Lambda(M_x', \sodot_{b_{iy} >0} M_i^{\snconv b_{iy}})-\Lambda(M_x',\sodot_{b_{iy} <0} M_i^{\snconv -b_{iy}})  \\
&&=\sum_{b_{iy} >0} \Lambda(M_x',M_i)b_{iy} -\sum_{b_{iy} <0} \Lambda(M_x',M_i)(-b_{iy})  \\
&&=\sum_{i\in\K} \Lambda(M_x', M_i)b_{iy}
=\sum_{i \neq x} \Lambda(M_x',M_i)b_{iy} + \Lambda(M_x',M_x)b_{xy} \\
&&=\sum_{i \neq x} \Lambda(M_x', M_i)(b_{iy}' -\delta(b_{ix} >0)b_{ix} b_{xy}) + \Lambda(M_x',M_x)b_{xy} \\
&&=\sum_{i \neq x} \Lambda(M_x', M_i) b_{iy}'
-\sum_{b_{ix} >0} \Lambda(M_x', M_i)b_{ix}b_{xy} + \Lambda(M_x',M_x)b_{xy} \\
&&\mathop=\limits_{\mathrm(a)}0-\Lambda(M_x', \sodot_{b_{ix} >0} M_i^{\snconv b_{ix}})b_{xy} + \Lambda(M_x',M_x)b_{xy} \\
&&=\bl-\Lambda(M_x',\sodot_{b_{ix} >0} M_i^{\snconv b_{ix}} ) + \Lambda(M_x',M_x) \br b_{xy} \\
&&=(-\Lambda(M_x', M_x' \hconv M_x) + \Lambda(M_x',M_x) )b_{xy} \\
&&=(-\Lambda(M_x', M_x')-\Lambda(M_x', M_x) + \Lambda(M_x',M_x) )b_{xy}=0.
\eneqn
Note that we used the compatibility of the pair 
$\bl\bl-\Lambda(\mu_x(M_i),\mu_x(M_j))\br_{i,j\in\K}, \, \mu_x(\widetilde B)\br$
when we derive the equality (a).

Since $L=(M'_x)^{\snconv b_{xy}}$, the equality $\La(M'_x,X)=\La(M'_x,Y)$ implies
\eqn
&& \Lambda(L,X)
=\Lambda(L,Y)=\La(L, M_y \conv M_y'). \label{eq:LaLX=LaLY}
\eneqn
Hence  the following  diagram is commutative:
\eq \label{eq:commutative ses}
&&\ba{c}\xymatrix{
0 \ar[r] & q L \conv X \ar[r] \ar[d]^{\rmat{L,X}} &  q^{\tLa(M_y,M_y')} L \conv (M_y \conv M_y')
 \ar[d]^{\rmat{L,M_y \circ M_y'}} \ar[r] &  L \conv Y \ar[r] \ar[d]_{\rmat{L,Y}} ^{\bwr}&0\\
0 \ar[r] & q^{d+1} X \conv L \ar[r]  & q^{d+\tLa(M_y,M_y')} (M_y \conv M_y') \conv L \ar[r] & q^{d}\, Y \conv L \ar[r] & 0,
}\ea
\eneq
where $d=-\Lambda(L,X)=-\Lambda(L,M_y \conv M_y')=-\Lambda(L,Y)$.
Note that  since $L=(M_x')^{\snconv b_{xy}}$ commutes with $Q$ and $A$, $\rmat{L,Y}$ is an isomorphism and hence
we have
$$\Im(\rmat{L,Y}) \simeq L \conv Y.$$
Hence we have an exact sequence
\eq
\xymatrix{0\ar[r]&\Im(\rmat{L,X})\ar[r]&\Im(\rmat{L, M_y\circ M'_y})\ar[r]
&L\circ Y\ar[r]&0.}\label{eq:exIm}
\eneq

On the other hand, $\rmat{L, M_y \circ M_y'}$ decomposes 
(up to a grading shift) as follows:
\eqn
\xymatrix@C=6em{
L \conv M_y \conv M_y'  \ar[r]^\sim_{\rmat{L,M_y} \circ M_y'} 
\ar@/^2pc/ [rr]^{\rmat{L, M_y \circ M_y'}}&
 M_y \conv L \conv M_y'  \ar[r]_{M_y \circ \rmat{L,M_y'}} &
 M_y  \conv M_y' \conv L.
 }
 \eneqn
Since $L=(M_x')^{\snconv b_{xy}}$ commutes with $M_y$,
the homomorphisms $\rmat{L,M_y} \conv M_y'$ is an isomorphism and hence we have
\eqn
\Im
 (\rmat{L, M_y \circ M_y'}) \simeq  M_y \conv (L\hconv M'_y)
\quad\text{up to a grading shift.}
 \eneqn
Similarly $\rmat{L, X}$ decomposes 
(up to a grading shift) as follows:
\eqn
\xymatrix@C=6em{
L \conv V \conv P  \ar[r]_{\rmat{L,V} \circ P} 
\ar@/^1.5pc/ [rr]^{\rmat{L,X}}&
 V \conv L \conv P  \ar[r]_{ V\circ \rmat{L,P}}^\sim &
V \conv P  \conv L.
 }
 \eneqn
Since $L$ commutes with $P$,
the homomorphisms  $V \conv \rmat{L,P} $ is an isomorphism and hence we have
  \eqn
\Im
 (\rmat{L,X})
  \simeq (L \hconv V) \conv P
\simeq \bl(M_x')^{\conv b_{xy}} \hconv M_x^{\conv b_{xy}}\br \conv P
\quad\text{up to a grading shift.}
 \eneqn
On the other hand, 
Lemma~\ref{lem:convpower} implies that
\eqn
(M_x')^{\conv b_{xy}} \hconv M_x^{\conv b_{xy}}
\simeq (M_x' \hconv M_x)^{\conv b_{xy}}
\simeq  C^{\conv b_{xy}} 
\simeq B \nconv A,
\eneqn
and hence we obtain
  \eqn
\Im
 (\rmat{L,X})
\simeq (B \nconv P)  \nconv A\quad\text{up to a grading shift.}
\eneqn
Thus the exact sequence \eqref{eq:exIm} becomes the exact sequence in $\shc$:
\eq
&&\xymatrix{
0 \ar[r] &  q^m (B \nconv P) \nconv A \ar[r]  &q^nM_y \conv (L\hconv M_y') \ar[r] & (L  \nconv Q) \nconv A \ar[r] & 0.
}\label{eq:exML}
\eneq
for some $m,n\in\Z$.
Since $(L  \nconv Q) \nconv A$ is self-dual, $n=\tL(M_y,L\hconv M_y')$.
On the other hand, we have
\eqn
\de(M_y,L \hconv M_y') \le \de(M_y, L) + \de(M_y, M_y') =1.
\eneqn
By the exact sequence \eqref{eq:exML},  
$M_y \conv (L\hconv M_y')$ is not simple and we conclude
$$\de(M_y,L \hconv M_y')=1.$$
Then Corollary~\ref{cor:crde} implies that $m=1$.
Thus we obtain an exact sequence in $\shc$:
\eq
&&\xymatrix@C=3ex{
0 \ar[r] &  q (B \nconv P) \nconv A \ar[r]  &q^{\tL(M_y,L\shconv M_y')}M_y \conv (L\hconv M_y') \ar[r] & (L  \nconv Q) \nconv A \ar[r] & 0.
}\label{eq:exML1}
\eneq
Now we shall rewrite \eqref{eq:exML1} by using $\scbul\conv A$ instead of 
$\scbul\nconv A$.
We have \eqn
&&\tLa(B,A)+\tLa(A,A)
=b_{xy}\tLa(C,A)
=b_{xy}\tLa(M_x' \hconv M_x,A) \\
&&\hs{10ex}=b_{xy}\tLa(M_x',A)+b_{xy}\tLa(M_x,A)=\tL(L,A)+b_{xy}\tLa(M_x,A).
\eneqn
On the other hand, the exact sequence \eqref{eq:sesMy} gives
\eqn
&&b_{xy}\tLa(M_x,A)+\tLa(P,A)=\tL(X,A)=\tL(M'_y\hconv M_y,A)\\
&&\hs{5ex}=\tL(M'_y,A)+\tL(M_y,A)=\tL(M_y\hconv M'_y,A)=
\tL(Y,A)=\tL(Q,A)+\tL(A,A).
\eneqn
It follows that
\eqn
&&\tLa(B \conv P,A)=\tLa(B,A) + \tLa(P,A) \\
&&\hs{5ex}=  \bl \tL(L,A)+b_{xy}\tLa(M_x,A)-\tLa(A,A)\br 
+\bl \tL(Q,A)+\tL(A,A)-b_{xy} \tLa(M_x,A)\br\\
&&\hs{5ex}=\tLa(L,A) + \tLa(Q,A)=\tLa(L \conv Q, A). \label{eq:qpower3}
\eneqn

Thus we have
\eq
&&\xymatrix{
0 \ar[r] &  q (B \nconv P) \conv A \ar[r]  &q^{c}  M_y \conv (L\hconv M_y') \ar[r] & (L  \nconv Q) \conv A \ar[r] & 0,
}\label{eq:LQA}
\eneq
where $c=\tLa(M_y, L \hconv M_y')-\tLa(B \nconv P,A)$.

Thus we obtain the identity in $K(R \gmod)$:
\eqn
q^c  [M_y] [L \hconv M_y'] = \bl q [B \nconv P] + [L \nconv Q]\br [A].
\eneqn

On the other hand,  hypothesis \eqref{eq:Cluster} implies that 
there exists $\phi \in \Q(q^{1/2}) \otimes_{\Z[q^{\pm1}]} K(\shc)$
such that
\eq
[M_y] \phi = q [B \nconv P] + [L \nconv Q]
\label{eq:Mphi}
\eneq
and
\eq
\phi [\mu_x(M)_i] = q^{\la'_{yi}} [\mu_x(M)_i] \phi
\quad
\text{for} \ i \neq y,
\label{eq:phiA}
\eneq
where $\mu_y\mu_x(-\La)=(\la'_{ij})_{i,j\in\K}$.

Hence, in $\Q(q^{1/2}) \otimes_{\Z[q^{\pm1}]} K(\shc)$, we have
\eqn
[M_y] \phi [A] = \bl q [B \nconv P] + [L \nconv Q]\br [A ]=  q^c  [M_y][L \hconv M_y'].
\eneqn
Since $\Q(q^{1/2}) \otimes_{\Z[q^{\pm1}]} K(\shc)$ is a domain, we conclude that
\eqn
\phi [A] = q^c [L \hconv M_y'].
\eneqn

On the other hand,  \eqref{eq:phiA} implies that
\eqn
\phi [A] = q^{l} [A] \phi\quad\text{for some $l\in\Z$.}
\eneqn
Hence, Theorem \ref{thm:divisible} implies that, when we write
\eqn
\phi =\sum_{b \in B(\infty)} a_b [L_b]\quad\text{for some $a_b \in \Q(q^{1/2})$,}
\eneqn
we have $L_b\conv A \simeq q^{l} A\conv L_b$ whenever $a_b\not=0$.
Hence we have
\eqn
q^c (L \hconv M_y')=\phi [A] = \sum_{b \in B(\infty)} a_b [L_b \conv A].
\eneqn
We conclude that
there exists  a self-dual simple module  $M_y''$ in $R \gmod$  such that
$M''_y$ commutes with $A$ and 
$$ \phi=q^m[M''_y]$$
for some $m\in\Z$.
Then \eqref{eq:Mphi} implies that
$$q^m[M_y\conv M''_y] = q [B \nconv P] + [L \nconv Q].$$
Hence there exists an exact sequence
$$0\To X\To q^m\;M_y\conv M''_y\To Y\To 0,$$
where $X=q B \nconv P$ and $Y=L \nconv Q$ or 
$X=L \nconv Q$ and $Y=q B \nconv P$.
By Corollary~\ref{cor:crde}, the last case does not occur and we have an exact sequence
$$0\To q B \nconv P\To q^m\;M_y\conv M''_y\To L \nconv Q\To 0.$$

Since $M_y$, $M_y''$ and $L \nconv Q$ are self-dual, we have
$m=\tLa(M_y,M_y'')$,
and we obtain the desired short exact sequence \eqref{eq:seqdes'}.

\vskip 1em
Since $\phi$ commutes with $[\mu_x (M)_i]$  up to a power of $q$ in
$K(\shc)$,
and
$\mu_x(M)_i$ is real simple,
$M_y''$ commutes with $\mu_x (M)_i$ for $i\neq y$.
\QED

\Cor\label{cor:main}
Let  $(\{M_i\}_{i\in\K},\widetilde B)$ be an admissible pair in $\shc$.
Under the assumption \eqref{eq:Cluster},
$\shc$ is a monoidal categorification of the quantum cluster algebra
$\mathscr A_{q^{1/2}}([\seed])$.
Furthermore, we have the followings{\rm:}
\bnum
  \item The quantum monoidal seed $\seed=(\{M_i\}_{i\in\K},  -\La,\widetilde B,D)$
  admits successive mutations in all directions.
  \item Any cluster monomial in $\Z[q^{\pm1/2}] \otimes_{\Z[q^\pm1]} K(\shc)$
   is the isomorphism class of a  real simple object in $\shc$ up to a power of $q^{1/2}$.
  \item Any cluster monomial in $\Z[q^{\pm1/2}]\otimes_{\Z[q^\pm1]} K(\shc)$
is a Laurent polynomial of the initial cluster variables with coefficient 
in $\Z_{\ge0}[q^{\pm1/2}]$.
\end{enumerate}
\encor

\Proof
(i) and (ii) are straightforward.

Let us show (iii). Let $x$ be a cluster monomial.
By the Laurent phenomenon (\cite{BZ05}),
we can
write 
$$xX^{\mathbf{c}}=\sum_{\bfa\in\Z_{\ge0}^\K}c_\bfa X^\bfa,$$
where 
$X=(X_i)_{i\in \K}$ is the initial cluster,  $\mathbf{c}\in\Z_{\ge0}^\K$
and $c_\bfa\in \Q(q^{\pm1/2})$.
Since $x$ and $X^{\mathbf{c}}$ are the isomorphism classes of simple modules
up to a power of $q^{1/2}$,
their product $xX^{\mathbf{c}}$ can be written as 
a linear combination of
the isomorphism classes of simple modules
with coefficients in $\Z_{\ge0}[q^{\pm1/2}]$.
Since
every $X^\bfa$ is the isomorphism class of a simple module
up to a power of $q^{1/2}$,
we have
$c_\bfa\in \Z_{\ge0}[q^{\pm1/2}]$.
\QED

\end{document}